\documentclass [12pt]{article}
\usepackage{amsmath,amssymb}

\begin{document}

   \title {Determinantal representations of the Moore-Penrose  inverse over
the quaternion skew field and corresponding Cramer's rules.}
\author {Ivan Kyrchei\footnote{Pidstrygach Institute
for Applied Problems of Mechanics and Mathematics, str.Naukova 3b,
Lviv, Ukraine, 79005, kyrchei@lms.lviv.ua}}
\date{}
\maketitle
\begin{abstract}
Determinantal representation of the Moore-Penrose inverse over the
quaternion skew field is obtained within the framework of a theory
of the column and row determinants. Using the obtained analogs of
the adjoint matrix, we get the Cramer rules for the least squares
solution of left and right systems of quaternionic linear
equations.
\end{abstract}
\textit{Keywords}: Moore-Penrose inverse, quaternion matrix, least
square solution, Cramer rule, quaternionic linear equation.

\noindent \textit{AMS classification}:15A06,15A09, 15A15, 15A33

\section{Introduction}
\newtheorem{definition}{Definition}[section]
\newtheorem{proposition}{Proposition}[section]
\newtheorem{lemma}{Lemma}[section]
\newtheorem{theorem}{Theorem}[section]
\newtheorem{corollary}{Corollary}[section]
\newtheorem{remark}{Remark}[section]
The existence, uniqueness and full-rank representation of the
Moore-Penrose inverse  over the quaternion skew field
${\rm{\mathbb{H}}}$ have been studied in particular in  \cite{ja,
wa, zho}.
 At the same time the problem of  determinantal representation
of the  Moore-Penrose inverse of a quaternion matrix is still
unsolved. Moreover, by means
 of any introduced determinants
currently it is not possible to receive a determinantal
representation of an inverse matrix over ${\rm{\mathbb{H}}}$. The
problem consists in the fact that there is no  determinant
functional over a skew field which would keep all the properties
inherent in the complex case, as is proved in \cite{as}. So the
determinants of Dieudonne \cite{di} and Study \cite{st} assume
values not in a skew field, but in the field which is its center.
These determinants can not be expanded by cofactors along an
arbitrary row or column of a matrix. The Moore determinant
 is introduced \cite{mo} in terms of permutations only on the class
of quaternion Hermitian matrices. The other determinant in terms
of permutations, the determinant of L. Chen \cite{ch} dissatisfies
a key property of a determinant, its singularity for noninvertible
matrices. The double determinant introduced by L. Chen can not be
expanded by cofactors along an arbitrary row or column of a matrix
as well.

In this paper determinantal representations of the Moore-Penrose
inverse are obtained within the framework of the theory of
  new  matrix functionals over the quaternion
skew field (the column and row determinants) introduced in
\cite{ky1}. In the first point we  cite some provisions from the
theory of the column and row determinants which are necessary for
the following. The theory of the column and row determinants of a
quaternionic matrix is considered completely  in \cite{ky1}.  In
the second point some known facts from the theory of eigenvalues
of a quaternion matrix and its singular value decomposition are
considered. The concept of a characteristic polynomial of a
Hermitian matrix is introduced and its coefficients are
investigated. A  Moore-Penrose  inverse  over a quaternion skew
field and its limit representation are introduced in this point as
well. In the third point the theorems about determinantal
representations of the Moore-Penrose inverse and the projection
matrices ${\rm {\bf A}}^{ +} {\rm {\bf A}}$, ${\rm {\bf A}}{\rm
{\bf A}}^{ +} $ for an arbitrary matrix ${\rm {\bf A}} \in {\rm
{\mathbb{H}}}^{m\times n}$  are proved. In the fourth point the
Cramer rules for a least squares solutions of  right and left
systems of linear equations over a quaternion skew field are
obtained.

\section{Elements of the theory of the column and row determinants.}

Let ${\rm M}\left( {n,{\rm {\mathbb{H}}}} \right)$ be the ring of
$n\times n$ quaternion matrices. By ${\rm {\mathbb{H}}}^{m\times
n}$ denote the set of all $m\times n$ matrices over the quaternion
skew field ${\rm {\mathbb{H}}}$ and by ${\rm
{\mathbb{H}}}_{r}^{m\times n} $ denote its subset of matrices of
rank $r$. Suppose $S_{n}$ is the symmetric group on the set
$I_{n}=\{1,\ldots,n\}$.
\begin{definition}
 The $i$th row determinant of ${\rm {\bf A}}=(a_{ij}) \in {\rm
M}\left( {n,{\mathbb{H}}} \right)$ is defined by

 \[{\rm{rdet}}_{ i} {\rm {\bf A}} =
{\sum\limits_{\sigma \in S_{n}} {\left( { - 1} \right)^{n -
r}{a_{i{\kern 1pt} i_{k_{1}}} } {a_{i_{k_{1}}   i_{k_{1} + 1}}}
\ldots } } {a_{i_{k_{1} + l_{1}}
 i}}  \ldots  {a_{i_{k_{r}}  i_{k_{r} + 1}}}
\ldots  {a_{i_{k_{r} + l_{r}}  i_{k_{r}} }}\] \noindent for all $i
= 1,\ldots,n $. The elements of the permutation $\sigma$ are
indices of each monomial. The left-ordered cycle notation of the
permutation $\sigma$ is written as follows,
\[\sigma = \left(
{i\,i_{k_{1}}  i_{k_{1} + 1} \ldots i_{k_{1} + l_{1}} }
\right)\left( {i_{k_{2}}  i_{k_{2} + 1} \ldots i_{k_{2} + l_{2}} }
\right)\ldots \left( {i_{k_{r}}  i_{k_{r} + 1} \ldots i_{k_{r} +
l_{r}} } \right).\] \noindent The index $i$ opens the first cycle
from the left  and other cycles satisfy the following conditions,
$i_{k_{2}} < i_{k_{3}}  < \ldots < i_{k_{r}}$ and $i_{k_{t}}  <
i_{k_{t} + s} $ for all $t = 2,\ldots,r $ and $s =1,\ldots,l_{t}
$.
\end{definition}

\begin{definition}
The $j$th column determinant
 of ${\rm {\bf
A}}=(a_{ij}) \in {\rm M}\left( {n,{\mathbb{H}}} \right)$ is
defined by
 \[{\rm{cdet}} _{{j}}\, {\rm {\bf A}} =
{{\sum\limits_{\tau \in S_{n}} {\left( { - 1} \right)^{n -
r}a_{j_{k_{r}} j_{k_{r} + l_{r}} } \ldots a_{j_{k_{r} + 1}
i_{k_{r}} }  \ldots } }a_{j\, j_{k_{1} + l_{1}} }  \ldots  a_{
j_{k_{1} + 1} j_{k_{1}} }a_{j_{k_{1}} j}}\] \noindent for all $j
=1,\ldots,n $. The right-ordered cycle notation of the permutation
$\tau \in S_{n}$ is written as follows,
 \[\tau =
\left( {j_{k_{r} + l_{r}}  \ldots j_{k_{r} + 1} j_{k_{r}} }
\right)\ldots \left( {j_{k_{2} + l_{2}}  \ldots j_{k_{2} + 1}
j_{k_{2}} } \right){\kern 1pt} \left( {j_{k_{1} + l_{1}}  \ldots
j_{k_{1} + 1} j_{k_{1} } j} \right).\] \noindent The index $j$
opens  the first cycle from the right  and other cycles satisfy
the following conditions, $j_{k_{2}}  < j_{k_{3}}  < \ldots <
j_{k_{r}} $ and $j_{k_{t}}  < j_{k_{t} + s} $ for all $t =
2,\ldots,r $ and $s = 1,\ldots,l_{t}  $.
\end{definition}

Suppose ${\rm {\bf A}}_{}^{i{\kern 1pt} j} $ denotes the submatrix
of ${\rm {\bf A}}$ obtained by deleting both the $i$th row and the
$j$th column. Let ${\rm {\bf a}}_{.j} $ be the $j$th column and
${\rm {\bf a}}_{i.} $ be the $i$th row of ${\rm {\bf A}}$. Suppose
${\rm {\bf A}}_{.j} \left( {{\rm {\bf b}}} \right)$ denotes the
matrix obtained from ${\rm {\bf A}}$ by replacing its $j$th column
with the column ${\rm {\bf b}}$, and ${\rm {\bf A}}_{i.} \left(
{{\rm {\bf b}}} \right)$ denotes the matrix obtained from ${\rm
{\bf A}}$ by replacing its $i$th row with the row ${\rm {\bf b}}$.

We  note some properties of column and row determinants of a
quaternion matrix ${\rm {\bf A}} = \left( {a_{ij}} \right)$, where
$i \in I_{n} $, $j \in J_{n} $ and $I_{n} = J_{n} = {\left\{
{1,\ldots ,n} \right\}}$.
\begin{proposition}  \cite{ky1}
If $b \in {\mathbb{H}}$, then
 $ {\rm{rdet}}_{ i} {\rm {\bf A}}_{i.} \left( {b
\cdot {\rm {\bf a}}_{i.}}  \right) = b \cdot {\rm{rdet}}_{ i} {\rm
{\bf A}}$ for all $i =1,\ldots,n $.
\end{proposition}
\begin{proposition} \cite{ky1}\label{kyrc8}
If $b \in {\mathbb{H}}$, then  ${\rm{cdet}} _{{j}}\, {\rm {\bf
A}}_{.j} \left( {{\rm {\bf a}}_{.j} b} \right) = {\rm{cdet}}
_{{j}}\, {\rm {\bf A}}  b$ for all $j =1,\ldots,n$.
\end{proposition}
\begin{proposition} \cite{ky1}
If for  ${\rm {\bf A}}\in {\rm M}\left( {n,{\mathbb{H}}}
\right)$\, there exists $t \in I_{n} $ such that $a_{tj} = b_{j} +
c_{j} $\, for all $j = 1,\ldots,n$, then
\[
\begin{array}{l}
   {\rm{rdet}}_{{i}}\, {\rm {\bf A}} = {\rm{rdet}}_{{i}}\, {\rm {\bf
A}}_{{t{\kern 1pt}.}} \left( {{\rm {\bf b}}} \right) +
{\rm{rdet}}_{{i}}\, {\rm {\bf A}}_{{t{\kern 1pt}.}} \left( {{\rm
{\bf c}}} \right), \\
  {\rm{cdet}} _{{i}}\, {\rm {\bf A}} = {\rm{cdet}} _{{i}}\, {\rm
{\bf A}}_{{t{\kern 1pt}.}} \left( {{\rm {\bf b}}} \right) +
{\rm{cdet}}_{{i}}\, {\rm {\bf A}}_{{t{\kern 1pt}.}} \left( {{\rm
{\bf c}}} \right),
\end{array}
\]
\noindent where ${\rm {\bf b}}=(b_{1},\ldots, b_{n})$, ${\rm {\bf
c}}=(c_{1},\ldots, c_{n})$ and for all $i =1,\ldots,n$.
\end{proposition}
\begin{proposition} \cite{ky1}
If for ${\rm {\bf A}}\in {\rm M}\left( {n,{\mathbb{H}}} \right)$\,
 there exists $t \in J_{n} $ such that $a_{i\,t} = b_{i} + c_{i}$
for all $i = 1,\ldots,n$, then
\[
\begin{array}{l}
  {\rm{rdet}}_{{j}}\, {\rm {\bf A}} = {\rm{rdet}}_{{j}}\, {\rm {\bf
A}}_{{\,.\,{\kern 1pt}t}} \left( {{\rm {\bf b}}} \right) +
{\rm{rdet}}_{{j}}\, {\rm {\bf A}}_{{\,.\,{\kern 1pt} t}} \left(
{{\rm
{\bf c}}} \right),\\
  {\rm{cdet}} _{{j}}\, {\rm {\bf A}} = {\rm{cdet}} _{{j}}\, {\rm
{\bf A}}_{{\,.\,{\kern 1pt}t}} \left( {{\rm {\bf b}}} \right) +
{\rm{cdet}} _{{j}} {\rm {\bf A}}_{{\,.\,{\kern 1pt}t}} \left(
{{\rm {\bf c}}} \right),
\end{array}
\]
\noindent where ${\rm {\bf b}}=(b_{1},\ldots, b_{n})^T$, ${\rm
{\bf c}}=(c_{1},\ldots, c_{n})^T$ and for all $j =1,\ldots,n$.
\end{proposition}

The following lemmas enable us to expand ${\rm{rdet}}_{{i}}\, {\rm
{\bf A}}$ by cofactors
  along  the $i$th row and ${\rm{cdet}} _{j} {\rm {\bf A}}$
 along  the $j$th column respectively for all $i, j = 1,\ldots,n$.

\begin{lemma}\label{kyrc1}  \cite{ky1}
Let $R_{i{\kern 1pt} j}$ be the right $ij$-th cofactor of ${\rm
{\bf A}}\in {\rm M}\left( {n, {\mathbb{H}}} \right)$, that is, $
{\rm{rdet}}_{{i}}\, {\rm {\bf A}} = {\sum\limits_{j = 1}^{n}
{{a_{i{\kern 1pt} j} \cdot R_{i{\kern 1pt} j} } }} $ for all $i =
1,\ldots,n$.  Then
\[
 R_{i{\kern 1pt} j} = {\left\{ {{\begin{array}{*{20}c}
  - {\rm{rdet}}_{{j}}\, {\rm {\bf A}}_{{.{\kern 1pt} j}}^{{i{\kern 1pt} i}} \left( {{\rm
{\bf a}}_{{.{\kern 1pt} {\kern 1pt} i}}}  \right),& {i \ne j},
\hfill \\
 {\rm{rdet}} _{{k}}\, {\rm {\bf A}}^{{i{\kern 1pt} i}},&{i = j},
\hfill \\
\end{array}} } \right.}
\]
\noindent where  ${\rm {\bf A}}_{.{\kern 1pt} j}^{i{\kern 1pt} i}
\left( {{\rm {\bf a}}_{.{\kern 1pt} {\kern 1pt} i}}  \right)$ is
obtained from ${\rm {\bf A}}$   by replacing the $j$th column with
the $i$th column, and then by deleting both the $i$th row and
column, $k = \min {\left\{ {I_{n}}  \right.} \setminus {\left.
{\{i\}} \right\}} $.
\end{lemma}
\begin{lemma} \cite{ky1}
Let $L_{i{\kern 1pt} j} $ be the left $ij$-th cofactor of
 ${\rm {\bf A}}\in {\rm M}\left( {n,{\mathbb{H}}} \right)$, that
 is,
$ {\rm{cdet}} _{{j}}\, {\rm {\bf A}} = {{\sum\limits_{i = 1}^{n}
{L_{i{\kern 1pt} j} \cdot a_{i{\kern 1pt} j}} }}$ for all $j
=1,\ldots,n$. Then
\[
L_{i{\kern 1pt} j} = {\left\{ {\begin{array}{*{20}c}
 -{\rm{cdet}} _{i}\, {\rm {\bf A}}_{i{\kern 1pt} .}^{j{\kern 1pt}j} \left( {{\rm {\bf a}}_{j{\kern 1pt}. } }\right),& {i \ne
j},\\
 {\rm{cdet}} _{k}\, {\rm {\bf A}}^{j\, j},& {i = j},
\\
\end{array} }\right.}
\]

\noindent where  ${\rm {\bf A}}_{i{\kern 1pt} .}^{jj} \left( {{\rm
{\bf a}}_{j{\kern 1pt} .} } \right)$ is obtained from ${\rm {\bf
A}}$
 by replacing the $i$th row with the $j$th row, and then by
deleting both the $j$th row and  column, $k = \min {\left\{
{J_{n}} \right.} \setminus {\left. {\{j\}} \right\}} $.
\end{lemma}
We recall some well-known definitions. The \textit{conjugate} of a
quaternion $a = a_{0} + a_{1} i + a_{2} j + a_{3} k \in {\rm
{\mathbb{H}}}$ is defined by $\overline {a} = a_{o}-a_{1}i- a_{2}j
- a_{3}k$. The \textit{Hermitian adjoint matrix} of ${\rm {\bf A}}
= \left( {a_{ij}} \right) \in {\rm {\mathbb{H}}}^{n\times m}$ is
called the matrix ${\rm {\bf A}}^{ *} = \left( {a_{ij}^{ *} }
\right)_{m\times n} $ if $a_{ij}^{ *} = \overline {a_{ji}}  $ for
all $i = 1,\ldots,n $ and $j = 1,\ldots,m$.
 The matrix ${\rm {\bf A}} = \left( {a_{ij}}  \right) \in {\rm
{\mathbb{H}}}^{n\times m}$ is \textit{Hermitian} if ${\rm {\bf
A}}^{ *}  = {\rm {\bf A}}$.

A following theorem has a key value in the theory of the column
and row determinants.

\begin{theorem} \cite{ky1}\label{kyrc2}
If ${\rm {\bf A}} = \left( {a_{ij}}  \right) \in {\rm M}\left(
{n,{\rm {\mathbb{H}}}} \right)$ is Hermitian, then ${\rm{rdet}}
_{1} {\rm {\bf A}} = \cdots = {\rm{rdet}} _{n} {\rm {\bf A}} =
{\rm{cdet}} _{1} {\rm {\bf A}} = \cdots = {\rm{cdet}} _{n} {\rm
{\bf A}} \in {\rm {\mathbb{R}}}.$
\end{theorem}
 Taking into account Theorem \ref{kyrc2} we define the determinant of a
Hermitian matrix by putting $\det {\rm {\bf A}}: = {\rm{rdet}}
_{i} {\rm {\bf A}} = {\rm{cdet}} _{i} {\rm {\bf A}}$ for all $i
=1,\ldots,n $. This determinant of a Hermitian matrix coincides
with the Moore determinant. The properties of the determinant of a
Hermitian matrix are considered in \cite{ky1} by means of the
column and row determinants. Among them we note the following.

\begin{theorem} \cite{ky1}
If the $i$th row of a Hermitian matrix ${\rm {\bf A}}\in {\rm
M}\left( {n,{\rm {\mathbb{H}}}} \right)$ is replaced with a left
linear combination of its other rows, i.e. ${\rm {\bf a}}_{i.} =
c_{1} {\rm {\bf a}}_{i_{1} .} + \cdots + c_{k}  {\rm {\bf
a}}_{i_{k} .}$, where $ c_{l} \in {\rm {\mathbb{H}}}$ for all $l
=1,\ldots,k$ and $\{i,i_{l}\}\subset I_{n} $, then

${\rm{cdet}}_{i} {\rm {\bf A}}_{i.} \left( {c_{1} \cdot {\rm {\bf
a}}_{i_{1} .} + \cdots + c_{k} \cdot {\rm {\bf a}}_{i_{k} .}}\,
\right) = {\rm{rdet}}_{i} {\rm {\bf A}}_{i.} \left( {c_{1} \cdot
{\rm {\bf a}}_{i_{1} .} + \cdots + c_{k} \cdot {\rm {\bf
a}}_{i_{k} .}}\, \right) = 0.$
\end{theorem}

\begin{theorem} \cite{ky1}
If the $j$th column of
 a Hermitian matrix ${\rm {\bf A}}\in
{\rm M}\left( {n,{\rm {\mathbb{H}}}} \right)$   is replaced with a
right linear combination of its other columns, i.e. ${\rm {\bf
a}}_{.j} = {\rm {\bf a}}_{.j_{1}}   c_{1} + \cdots + {\rm {\bf
a}}_{.j_{k}} c_{k} $, where $c_{l} \in {\rm {\mathbb{H}}}$ for all
$l =1,\ldots,k$ and $\{j,j_{l}\}\subset J_{n}$,  then

${\rm{cdet}}_{i} {\rm {\bf A}}_{.\,i} \left( {{\rm {\bf
a}}_{.\,i_{1}} \cdot c_{1} + \cdots + {\rm {\bf a}}_{.\,i_{k}}
\cdot c_{k}} \right) = {\rm{rdet}}_{i} {\rm {\bf A}}_{.\,i} \left(
{{\rm {\bf a}}_{.\,i_{1}} \cdot c_{1} + \cdots + {\rm {\bf
a}}_{.\,i_{k}} \cdot c_{k} } \right) = 0.$
\end{theorem}
The following theorem about determinantal representation of an
inverse matrix of Hermitian follows immediately from these
properties.

\begin{theorem}
 \cite{ky1}\label{kyrc10}
 If a Hermitian matrix ${\rm
{\bf A}} \in {\rm M}\left( {n,{\rm {\mathbb{H}}}} \right)$ is such
that $\det {\rm {\bf A}} \ne 0$, then there exist a unique right
inverse  matrix $(R{\rm {\bf A}})^{ - 1}$ and a unique left
inverse matrix $(L{\rm {\bf A}})^{ - 1}$, and $\left( {R{\rm {\bf
A}}} \right)^{ - 1} = \left( {L{\rm {\bf A}}} \right)^{ - 1} =
:{\rm {\bf A}}^{ - 1}$. They possess the following determinantal
representations:
\[
  \left( {R{\rm {\bf A}}} \right)^{ - 1} = {\frac{{1}}{{\det {\rm
{\bf A}}}}}
\begin{pmatrix}
  R_{11} & R_{21} & \cdots & R_{n1} \\
  R_{12} & R_{22} & \cdots & R_{n2} \\
  \cdots & \cdots & \cdots & \cdots \\
  R_{1n} & R_{2n} & \cdots & R_{nn}
\end{pmatrix},
  \left( {L{\rm {\bf A}}} \right)^{ - 1} = {\frac{{1}}{{\det {\rm
{\bf A}}}}}
\begin{pmatrix}
  L_{11} & L_{21} & \cdots  & L_{n1} \\
  L_{12} & L_{22} & \cdots  & L_{n2} \\
  \cdots  & \cdots  & \cdots  & \cdots  \\
  L_{1n} & L_{2n} & \cdots  & L_{nn}
\end{pmatrix}.
\]
Here $R_{ij}$,  $L_{ij}$ are right and left $ij$-th cofactors of
${\rm {\bf
 A}}$ respectively for all $i,j =1,\ldots,n$.
\end{theorem}
 To obtain
determinantal representation of an arbitrary inverse matrix ${\rm
{\bf A}}^{ - 1}$, we consider the right ${\rm {\bf A}}{\rm {\bf
A}}^{ *} $ and left ${\rm {\bf A}}^{ *} {\rm {\bf A}}$
corresponding Hermitian matrix.

\begin{theorem} \cite{ky1}\label{kyrc12}
If an arbitrary column  of ${\rm {\bf A}}\in {\rm
{\mathbb{H}}}^{m\times n} $ is a right  linear combination of its
other columns, or an arbitrary row of ${\rm {\bf A}}^{ * }$ is a
left linear combination of its others, then $\det {\rm {\bf A}}^{
* }{\rm {\bf A}} = 0.$
\end{theorem}
Since the principal submatrices of a Hermitian matrix  are
Hermitian, the principal minor may be defined as the determinant
of its principal submatrix by analogy to the commutative case. We
introduce \textit{the rank by principal minors} that is the
maximal order of a nonzero principal minor of a Hermitian matrix.
The following theorem determines a relationship between it and the
rank of a matrix defining as ceiling amount of right-linearly
independent columns  or left-linearly independent rows which form
basis.

\begin{theorem} \cite{ky1}\label{kyrc6}
A rank by principal minors of ${\rm {\bf A}}^{ *} {\rm {\bf A}}$
is equal to its rank and a rank of ${\rm {\bf A}} \in {\rm
{\mathbb{H}}}^{m\times n}$.
\end{theorem}
\begin{theorem} \cite{ky1}\label{kyrc11}
If ${\rm {\bf A}} \in {\rm {\mathbb{H}}}^{m\times n}$, then an
arbitrary column of  ${\rm {\bf A}}$ is a right linear combination
of its basis columns or an arbitrary row of  ${\rm {\bf A}}$ is a
left linear combination of its basis rows.
\end{theorem}
The criterion of singularity of a Hermitian matrix is obtained.

\begin{theorem} \cite{ky1}\label{kyrc13}
The right-linearly independence of columns of ${\rm {\bf A}} \in
{\rm {\mathbb{H}}}^{m\times n}$ or the left-linearly independence
of rows of ${\rm {\bf A}}^{ *} $ is the necessary and sufficient
condition for $\det {\rm {\bf A}}^{ *} {\rm {\bf A}} \neq 0.$
\end{theorem}
\begin{theorem}\cite{ky1}
If ${\rm {\bf A}} \in {\rm M}\left( {n,{\rm {\mathbb{H}}}}
\right)$, then $\det {\rm {\bf A}}{\rm {\bf A}}^{ *} = \det {\rm
{\bf A}}^{ *} {\rm {\bf A}}$.
\end{theorem}

A concept of the double determinant is introduced by this theorem.
This concept was initially introduced by L. Chen in \cite{ch}.

\begin{definition}
The determinant of the corresponding Hermitian matrix of ${\rm
{\bf A}} \in {\rm M}\left( {n,{\rm {\mathbb{H}}}} \right)$
 is called its double determinant, i.e.
${\rm{ddet}}{ \rm{\bf A}}: = \det \left( {{\rm {\bf A}}^{ *} {\rm
{\bf A}}} \right) = \det \left( {{\rm {\bf A}}{\rm {\bf A}}^{ *} }
\right)$.
\end{definition}

The relationship between the double determinant and the
noncommutative determinants of E. Moore, E. Study and J. Diedonne
is obtained, ${\rm ddet} {\rm {\bf A}} = {\rm Mdet} \left( {{\rm
{\bf A}}^{ *} {\rm {\bf A}}} \right) = {\rm Sdet} {\rm {\bf A}} =
{\rm Ddet} ^{2}{\rm {\bf A}}$. But unlike those, the double
determinant can be expanded  along an arbitrary row or column  by
means of the column and row determinants.

\begin{definition}
Suppose ${\rm {\bf A}} \in {\rm M}\left( {n,{\rm {\mathbb{H}}}}
\right)$. We have a column expansion of ${\rm{ddet}} {\rm {\bf
A}}$ along the $j$th column, ${\rm{ddet}} {\rm {\bf A}} =
{\rm{cdet}} _{j} \left( {{\rm {\bf A}}^{ *} {\rm {\bf A}}} \right)
= {\sum\limits_{i} {{\mathbb{L}} _{ij} \cdot a_{ij}} }$, and a row
expansion of it along the $i$th row,
 ${\rm{ddet}} {\rm {\bf A}} = {\rm{rdet}}_{i} \left(
{{\rm {\bf A}}{\rm {\bf A}}^{ *} } \right) = {\sum\limits_{j}
{a_{ij} \cdot} } {\mathbb{R}} _{ i{\kern 1pt}j} $ for all $i,j
=1,\ldots,n $. Then by definition of the left double $ij$th
cofactor we put ${\mathbb{L}} _{ij} $ and by definition of the
right double $ij$th cofactor we put ${\mathbb{R}} _{ i{\kern1pt}j}
$.
\end{definition}

\begin{theorem} \cite{ky1}
The necessary and sufficient condition of invertibility of  ${\rm
{\bf A}} \in {\rm M}(n,{\rm {\mathbb{H}}})$ is ${\rm{ddet}} {\rm
{\bf A}} \ne 0$. Then there exists $ {\rm {\bf A}}^{ - 1} = \left(
{L{\rm {\bf A}}} \right)^{ - 1} = \left( {R{\rm {\bf A}}}
\right)^{ - 1}$, where
\begin{equation}
\label{kyr1} \left( {L{\rm {\bf A}}} \right)^{ - 1} =\left( {{\rm
{\bf A}}^{ *}{\rm {\bf A}} } \right)^{ - 1}{\rm {\bf A}}^{ *}
={\frac{{1}}{{{\rm{ddet}}{ \rm{\bf A}} }}}
\begin{pmatrix}
  {\mathbb{L}} _{11} & {\mathbb{L}} _{21}& \ldots & {\mathbb{L}} _{n1} \\
  {\mathbb{L}} _{12} & {\mathbb{L}} _{22} & \ldots & {\mathbb{L}} _{n2} \\
  \ldots & \ldots & \ldots & \ldots \\
 {\mathbb{L}} _{1n} & {\mathbb{L}} _{2n} & \ldots & {\mathbb{L}} _{nn}
\end{pmatrix},
\end{equation}
\begin{equation}\label{kyr2} \left( {R{\rm {\bf A}}} \right)^{ - 1} = {\rm {\bf
A}}^{ *} \left( {{\rm {\bf A}}{\rm {\bf A}}^{ *} } \right)^{ - 1}
= {\frac{{1}}{{{\rm{ddet}}{ \rm{\bf A}} }}}
\begin{pmatrix}
 {\mathbb{R}} _{\,{\kern 1pt} 11} & {\mathbb{R}} _{\,{\kern 1pt} 21} &\ldots & {\mathbb{R}} _{\,{\kern 1pt} n1} \\
 {\mathbb{R}} _{\,{\kern 1pt} 12} & {\mathbb{R}} _{\,{\kern 1pt} 22} &\ldots & {\mathbb{R}} _{\,{\kern 1pt} n2}  \\
 \ldots  & \ldots & \ldots & \ldots \\
 {\mathbb{R}} _{\,{\kern 1pt} 1n} & {\mathbb{R}} _{\,{\kern 1pt} 2n} &\ldots & {\mathbb{R}} _{\,{\kern 1pt} nn}
\end{pmatrix},
\end{equation}

and ${\mathbb{L}} _{ij} = {\rm{cdet}} _{j} ({\rm {\bf
A}}^{\ast}{\rm {\bf A}})_{.j} \left( {{\rm {\bf a}}_{.{\kern 1pt}
i}^{ *} } \right)$, ${\mathbb{R}} _{\,{\kern 1pt} ij} =
{\rm{rdet}}_{i} ({\rm {\bf A}}{\rm {\bf A}}^{\ast})_{i.} \left(
{{\rm {\bf a}}_{j.}^{ *} }  \right)$  for all $i,j =1,\ldots,n$.
\end{theorem}

This theorem introduces the determinantal representations of an
inverse matrix by the left (\ref{kyr1}) and right (\ref{kyr2})
double cofactors.

\section{ The singular value decomposition
 and the Moore-Penrose inverse of a quaternion matrix}

Due to the noncommutativity of quaternions, there are two types of
eigenvalues.
\begin{definition}
Let ${\rm {\bf A}} \in {\rm M}\left( {n,{\rm {\mathbb{H}}}}
\right)$. A quaternion $\lambda$ is said to be a right eigenvalue
of ${\rm {\bf A}}$ if ${\rm {\bf A}} \cdot {\rm {\bf x}} = {\rm
{\bf x}} \cdot \lambda $ for some nonzero quaternion column-vector
${\rm {\bf x}}$. Similarly $\lambda$ is a left eigenvalue if ${\rm
{\bf A}} \cdot {\rm {\bf x}} = \lambda \cdot {\rm {\bf x}}$.
\end{definition}
The theory on the left eigenvalues of quaternion matrices has been
investigated in particular in \cite{hu, so, wo}. The theory on the
right eigenvalues of quaternion matrices is more developed. In
particular we note  \cite{ba, dr, zha}. From this theory we cite
the following propositions.

\begin{proposition}\cite{zha}
Let ${\rm {\bf A}} \in {\rm M}\left( {n,{\rm {\mathbb{H}}}}
\right)$ is Hermitian. Then ${\rm {\bf A}}$ has exactly $n$ real
right eigenvalues.
\end{proposition}
\begin{definition}
Suppose ${\rm {\bf U}} \in {\rm M}\left( {n,{\rm {\mathbb{H}}}}
\right)$. If ${\rm {\bf U}}^{ *} {\rm {\bf U}} = {\rm {\bf U}}{\rm
{\bf U}}^{ *}  = {\rm {\bf I}}$, then the matrix ${\rm {\bf U}}$
is called unitary.
\end{definition}

\begin{proposition} \cite{zha}
Let ${\rm {\bf A}} \in {\rm M}\left( {n,{\rm {\mathbb{H}}}}
\right)$ be given. Then, ${\rm {\bf A}}$ is Hermitian  if and only
if there are a unitary matrix ${\rm {\bf U}} \in {\rm M}\left(
{n,{\rm {\mathbb{H}}}} \right)$ and a real diagonal matrix ${\rm
{\bf D}} = {\rm diag}\left( {\lambda _{{\kern 1pt} 1} ,\lambda
_{{\kern 1pt} 2} ,\ldots ,\lambda _{{\kern 1pt} n}}  \right)$ such
that ${\rm {\bf A}} = {\rm {\bf U}}{\rm {\bf D}}{\rm {\bf U}}^{
*}$, where $\lambda _{ 1},...,\lambda _{ n} $ are right
eigenvalues of ${\rm {\bf A}}$.
\end{proposition}
Suppose ${\rm {\bf A}} \in {\rm M}\left( {n,{\rm {\mathbb{H}}}}
\right) $ is Hermitian and $\lambda \in {\rm {\mathbb {R}}}$ is a
right eigenvalue, then ${\rm {\bf A}} \cdot {\rm {\bf x}} = {\rm
{\bf x}} \cdot \lambda = \lambda \cdot {\rm {\bf x}}$. This means
that all right eigenvalues of a Hermitian matrix are its left
eigenvalues as well. For real left eigenvalues, $\lambda \in {\rm
{\mathbb {R}}}$, the matrix $\lambda {\rm {\bf I}} - {\rm {\bf
A}}$ is Hermitian.
\begin{definition}
If $t \in {\rm {\mathbb {R}}}$, then for a Hermitian matrix ${\rm
{\bf A}}$ the polynomial $p_{{\rm {\bf A}}}\left( {t} \right) =
\det \left( {t{\rm {\bf I}} - {\rm {\bf A}}} \right)$ is said to
be the characteristic polynomial of ${\rm {\bf A}}$.
\end{definition}
The roots of the characteristic polynomial of a Hermitian matrix
are its real left eigenvalues, which are its right eigenvalues as
well. We shall investigate coefficients of the characteristic
polynomial as in the commutative case (see, e.g. \cite{la}). At
first we prove the auxiliary lemma.
\begin{lemma}\label{kyrc3}
Let ${\rm {\bf A}} \in {\rm M}\left( {n,{\rm \mathbb{H}}} \right)$
be Hermitian and the columns $i_{1} ,\ldots ,i_{k} $ of ${\rm {\bf
A}}$ coincide with the unit vectors ${\rm {\bf e}}_{i_{1}} ,\ldots
,{\rm {\bf e}}_{i_{k}}$. Then $\det {\rm {\bf A}}$  equals  a
principal minor obtained from ${\rm {\bf A}}$ by deleting the rows
and columns $i_{1} ,\ldots ,i_{k}$.

\end{lemma}
{\textit{Proof}}. We claim that if ${\rm {\bf A}} \in {\rm
M}\left( {n,{\rm {\mathbb{H}}}} \right)$ is Hermitian and the
columns $i_{1} ,\ldots ,i_{k} $ of ${\rm {\bf A}}$ coincide with
the unit column vectors ${\rm {\bf e}}_{i_{1}} ,\ldots ,{\rm {\bf
e}}_{i_{k}}$ respectively, then the rows $i_{1} ,\ldots ,i_{k} $
coincide with the unit row vectors ${\rm {\bf e}}_{i_{1}} ,\ldots
,{\rm {\bf e}}_{i_{k}}$ as well. Using Lemma \ref{kyrc1} we expand
$\det {\rm {\bf A}}$ along the  $i_{{\kern 1pt} 1} $th column,
where $a_{i_{1} {\kern 1pt} k} = 0$ for all $k \ne i_{1} $ and
$a_{i_{1} {\kern 1pt} i_{1}}  = 1$.
 Then we obtain
\[\begin{array}{c}
  \det {\rm {\bf A}} ={\rm{cdet}} _{i_{1}}  {\rm {\bf A}}=\\
   = - {\rm{cdet}} _{i_{1}}  {\rm {\bf A}}_{i_{1} {\kern 1pt} .}^{11} \left( {{\rm
{\bf a}}_{1{\kern 1pt} .}}  \right) \cdot 0 + \cdots + {\rm{cdet}}
_{1} {\rm {\bf A}}^{i_{1} i_{1}}  \cdot 1 + \cdots - {\rm{cdet}}
_{i_{1}}  {\rm {\bf A}}_{i_{1} {\kern 1pt} .}^{n{\kern 1pt} n}
\left( {{\rm {\bf a}}_{n{\kern 1pt} .}} \right) \cdot 0 =\\=
{\rm{cdet}} _{1} {\rm {\bf A}}^{i_{1} i_{1}}.
\end{array}
\]
Since the submatrix ${\rm {\bf A}}^{i_{1} i_{1}} $ is obtained
from ${\rm {\bf A}}$ by deleting both the $i_{1} $-th rows and
columns, by Theorem \ref{kyrc2} it follows that ${\rm{cdet}} _{1}
{\rm {\bf A}}^{i_{1} i_{1}} = \det {\rm {\bf A}}^{i_{1} i_{1}} $.
Now we calculate this principal minor expanding along the
$i_{{\kern 1pt} 2} $-th column. Similarly to above we have that
$\det {\rm {\bf A}}$ is equal to a principal minor obtained from
${\rm {\bf A}}$ by deleting both the $i_{1}$th and $i_{2}$th rows
and columns. Continuing this line of reasoning we complete the
proof of the lemma.$\blacksquare$

Taking into account Lemma \ref{kyrc3} we can prove the following
theorem by analogy to the commutative case (see, e.g. \cite{la}).

\begin{theorem}\label{kyrc7}
If ${\rm {\bf A}} \in {\rm M}\left( {n,{\rm {\mathbb{H}}}}
\right)$ is Hermitian, then $p_{{\rm {\bf A}}}\left( {t} \right) =
t^{n} - d_{1} t^{n - 1} + d_{2} t^{n - 2} - \cdots + \left( { - 1}
\right)^{n}d_{n}$, where $d_{r} $ is the sum of principle minors
of ${\rm {\bf A}}$ of order $r$, $1 \le r < n$, and $d_{n}=\det
{\rm {\bf A}}$.
\end{theorem}
For any quaternion matrix ${\rm {\bf A}} \in {\rm M}\left( {n,{\rm
{\mathbb{H}}}} \right)$,  the eigenvalues of  ${\rm {\bf A}}^{ *}
{\rm {\bf A}}$ are all nonnegative real numbers \cite{wi}.

\begin{definition}
Let ${\rm {\bf A}} \in {\rm {\mathbb{H}}}^{m\times n}$. The
nonnegative square roots of the $n$ eigenvalues of  ${\rm {\bf
A}}^{ *} {\rm {\bf A}}$ are called the singular values of ${\rm
{\bf A}}$.
\end{definition}
A key value for a determinantal representation of the
Moore-Penrose inverse over the quaternion skew field is the
following singular value decomposition (SVD) theorem.

\begin{theorem} \cite{wi, zha}
(SVD) Let ${\rm {\bf A}} \in {\rm {\mathbb{H}}}_{r}^{m\times n} $.
Then there exist unitary quaternion matrices ${\rm {\bf U}}_{1}
\in {\rm {\mathbb{H}}}^{m\times m}$ and ${\rm {\bf U}}_{2} \in
{\rm {\mathbb{H}}}^{n\times n}$ such that

\begin{equation}
\label{kyr3} {\rm {\bf U}}_{1} {\rm {\bf A}}{\rm {\bf U}}_{2} =
{\left[ {{\begin{array}{*{20}c}
 {{\rm {\bf D}}_{r}}  \hfill & {{\rm {\bf 0}}} \hfill \\
 {{\rm {\bf 0}}} \hfill & {{\rm {\bf 0}}} \hfill \\
\end{array}} } \right]} \in {\rm {\mathbb{H}}}^{m\times n},
\end{equation}
where ${\rm {\bf D}}_{r} = {\rm diag}\left( {\sigma _{1} ,\sigma
_{2} ,\ldots ,\sigma _{r}}  \right), \sigma _{1} \ge \sigma _{2}
\ge \cdots \ge \sigma _{r} > 0$, and $\sigma _{1} ,\sigma _{2}
,\ldots ,\sigma _{r} $ are the all nonzero singular values of
${\rm {\bf A}}.$
\end{theorem}
As unitary matrices are invertible, the equality (\ref{kyr3})
 can be written as follows
\begin{equation}
\label{kyr18}
 {\rm {\bf A}} = {\rm {\bf V}}{\rm {\bf \Sigma}
}{\rm {\bf W}}^{ *},
\end{equation}
 where ${\rm {\bf V}} \in {\rm
{\mathbb{H}}}^{m\times m}$ and ${\rm {\bf W}} \in {\rm
{\mathbb{H}}}^{n\times n}$ are unitary matrices, and the matrix
${\rm {\bf \Sigma} } = \left( {\sigma _{ij}} \right) \in {\rm
{\mathbb{H}}}_{r}^{m\times n} $ is such that $\sigma _{11} \ge
\sigma _{22} \ge \cdots \ge \sigma _{rr} > \sigma _{r + 1\,r + 1}
= \cdots = \sigma _{qq} = 0$, $q = \min {\left\{ {n,m} \right\}}$.

We get the following lemmas, which have the analogues in the
complex case \cite{ca}.

\begin{lemma}\label{kyrc4}
Suppose ${\rm {\bf A}} \in {\rm {\mathbb{H}}}^{m\times n}$ has the
singular value decomposition, ${\rm {\bf A}} = {\rm {\bf V}}{\rm
{\bf \Sigma} }{\rm {\bf W}}^{ *} $. Let ${\rm {\bf A}}^{ + } =
{\rm {\bf W}} \cdot {\rm {\bf \Sigma} }^{ +} \cdot {\rm {\bf V}}^{
* }$, where ${\rm {\bf \Sigma} }^{ +}  \in {\rm
{\mathbb{H}}}^{n\times m}$ is obtained from ${\rm {\bf \Sigma }}$
by transposition and replacing positive entries of ${\rm {\bf
\Sigma }}$ by reciprocal. Then for ${\rm {\bf A}}^{ +} $ the
following conditions are true
\begin{equation}\label{kyr16}
\begin{array}{l}
  1)\,\, \left( {{\rm {\bf A}}{\rm {\bf
A}}^{ +} } \right)^{ *}  = {\rm {\bf A}}{\rm {\bf A}}^{ +};\\
  2)\,\, \left( {{\rm {\bf A}}^{ +} {\rm {\bf A}}} \right)^{ *}  = {\rm
{\bf A}}^{ +} {\rm {\bf A}};\\
  3)\,\, {\rm {\bf A}}{\rm {\bf A}}^{ +}
{\rm {\bf A}} = {\rm {\bf A}};\\
  4)\,\,{\rm {\bf A}}^{ +} {\rm {\bf
A}}{\rm {\bf A}}^{ +}  = {\rm {\bf A}}^{ +}.
\end{array}
\end{equation}
\end{lemma}

{\textit{Proof}}. We obviously have $\left( {{\rm {\bf \Sigma}
}^{T}} \right)^{ *}  = {\rm {\bf \Sigma} }$ and $\left( {\left(
{{\rm {\bf \Sigma} }^{ +} } \right)^{T}} \right)^{ *}  = {\rm {\bf
\Sigma} }^{ +} $ for ${\rm {\bf \Sigma} }$ from the SVD by
(\ref{kyr18}) and ${\rm {\bf \Sigma} }^{ +} $. Then it follows
that

\[\begin{array}{c}
 \left( {{\rm {\bf A}}{\rm {\bf A}}^{ +} } \right)^{ *}  = \left(
{{\rm {\bf V}}{\rm {\bf \Sigma} }  {\rm {\bf W}}^{ *} {\rm {\bf
W}}  {\rm {\bf \Sigma} }^{ +} {\rm {\bf V}}^{ *} } \right)^{ *}  =
\left( {{\rm {\bf V}}{\rm {\bf \Sigma} } {\rm {\bf I}}  {\rm {\bf
\Sigma} }^{ + }{\rm {\bf V}}^{ *} } \right)^{ *}  =  \left( {{\rm
{\bf V}}  \left( {{\rm {\bf \Sigma} }^{ +} } \right)^{T}{\rm {\bf
\Sigma} }^{T}  {\rm
{\bf V}}^{ *} } \right)^{ *}  =\\
 =\left( {{\rm {\bf V}}\left( {{\rm
{\bf \Sigma} }^{ +} } \right)^{T}{\rm {\bf W}}^{ *} {\rm {\bf
W}}{\rm {\bf \Sigma} }^{T}{\rm {\bf V}}^{ *} } \right)^{ *} = {\rm
{\bf V}}  {\rm {\bf \Sigma} }{\rm {\bf W}}^{ *} {\rm {\bf W}} {\rm
{\bf \Sigma} }^{ +} {\rm {\bf V}}^{ *}  = {\rm {\bf A}}{\rm {\bf
A}}^{ +}.
\end{array}
\]
The proof of 1) is completed. By analogy we can prove 2).

Now we prove the condition 3). Note that ${\rm {\bf \Sigma} }{\rm
{\bf \Sigma} }^{ +}  = {\left[ {{\begin{array}{*{20}c}
 {{\rm {\bf I}}_{r}}  \hfill & {{\rm {\bf 0}}} \hfill \\
 {{\rm {\bf 0}}} \hfill & {{\rm {\bf 0}}} \hfill \\
\end{array}} } \right]} \in {\rm {\mathbb{H}}}^{m\times m}$. This implies ${\rm {\bf
\Sigma} }{\rm {\bf \Sigma} }^{ +} {\rm {\bf \Sigma} } = {\rm {\bf
\Sigma }}$, then $ {\rm {\bf A}}{\rm {\bf A}}^{ +} {\rm {\bf A}} =
{\rm {\bf V}}{\rm {\bf \Sigma} }{\rm {\bf W}}^{ *}  \cdot {\rm
{\bf W}}{\rm {\bf \Sigma} }^{ + }{\rm {\bf V}}^{ *}  \cdot {\rm
{\bf V}}{\rm {\bf \Sigma} }{\rm {\bf W}}^{ * } = {\rm {\bf V}}
\cdot {\rm {\bf \Sigma} }{\rm {\bf \Sigma} }^{ +} {\rm {\bf
\Sigma} } \cdot {\rm {\bf W}}^{ *}  = {\rm {\bf V}} \cdot {\rm
{\bf \Sigma} } \cdot {\rm {\bf W}}^{ *}  = {\rm {\bf A}}. $

By analogy to 3) can be prove the condition 4).$\blacksquare$

\begin{lemma}
There exists a unique matrix ${\rm {\bf A}}^{ +} $ that satisfies
conditions 1)-4) in  (\ref{kyr16}).
\end{lemma}
{\textit{Proof}}. Suppose that both matrices ${\rm {\bf B}} \in
{\rm {\mathbb{H}}}^{n\times m}$ and ${\rm {\bf C}} \in {\rm
{\mathbb{H}}}^{n\times m}$ satisfy conditions 1)-4) in Lemma
\ref{kyrc4}. Then we have
\[
\begin{array}{c}
   {\rm {\bf B}} = {\rm {\bf B}}{\rm {\bf A}}{\rm {\bf B}} = {\rm {\bf
B}}\left( {{\rm {\bf A}}{\rm {\bf B}}} \right)^{ *}  = {\rm {\bf
B}}{\rm {\bf B}}^{ *} {\rm {\bf A}}^{ *}  = {\rm {\bf B}}{\rm {\bf
B}}^{ *} \left( {{\rm {\bf A}}{\rm {\bf C}}{\rm {\bf A}}}
\right)^{ *}  = {\rm {\bf B}}{\rm {\bf B}}^{ *} {\rm {\bf A}}^{ *}
{\rm {\bf C}}^{ *} {\rm {\bf A}}^{ *}  =  \\
 ={\rm {\bf B}}\left( {{\rm
{\bf A}}{\rm {\bf B}}} \right)^{ *} \left( {{\rm {\bf A}}{\rm {\bf
C}}} \right)^{ *}  = {\rm {\bf B}}{\rm {\bf A}}{\rm {\bf B}} {\rm
{\bf A}}{\rm {\bf C}} = {\rm {\bf B}}{\rm {\bf A}}{\rm {\bf C}}  =
{\rm {\bf B}}{\rm {\bf A}}  {\rm {\bf C}}{\rm {\bf A}}{\rm {\bf
C}} = \left( {{\rm {\bf B}}{\rm {\bf A}}} \right)^{ *}  \left(
{{\rm {\bf C}}{\rm {\bf A}}} \right)^{ *} {\rm {\bf C}} =  \\
 ={\rm
{\bf A}}^{ *} {\rm {\bf B}}^{ *} {\rm {\bf A}}^{ *} {\rm {\bf
C}}^{ *} {\rm {\bf C}} =\left( {{\rm {\bf A}}{\rm {\bf B}}{\rm
{\bf A}}} \right)^{ *} {\rm {\bf C}}^{ *} {\rm {\bf C}} = {\rm
{\bf A}}^{ *} {\rm {\bf C}}^{ *} {\rm {\bf C}} = \left( {{\rm {\bf
C}}{\rm {\bf A}}} \right)^{ *} {\rm {\bf C}} = {\rm {\bf C}}{\rm
{\bf A}}{\rm {\bf C}} = {\rm {\bf C}}.
\end{array}
\]
$\blacksquare$

\begin{definition}
Let ${\rm {\bf A}} \in {\rm {\mathbb{H}}}^{m\times n}$. The matrix
${\rm {\bf A}}^{ +}$ is called the Moore-Penrose  inverse if it
satisfies all conditions  in (\ref{kyr16}).

\end{definition}

By analogy to the complex case \cite{ca} we have the theorem about
the limit representation of the Moore-Penrose inverse.

\begin{theorem}\label{kyrc9}
If ${\rm {\bf A}} \in {\rm {\mathbb{H}}}^{m\times n}$ and ${\rm
{\bf A}}^{ +}$ is its Moore-Penrose  inverse, then ${\rm {\bf
A}}^{ +}  = {\mathop {\lim} \limits_{\alpha \to 0}} {\rm {\bf
A}}^{ * }\left( {{\rm {\bf A}}{\rm {\bf A}}^{ *}  + \alpha {\rm
{\bf I}}} \right)^{ - 1} = {\mathop {\lim} \limits_{\alpha \to 0}}
\left( {{\rm {\bf A}}^{
* }{\rm {\bf A}} + \alpha {\rm {\bf I}}} \right)^{ - 1}{\rm {\bf
A}}^{ *} $, where $\alpha \in {\rm {\mathbb {R}}}_{ +}  $.
\end{theorem}
{\textit{Proof}}. Suppose ${\rm {\bf A}} = {\rm {\bf V}} {\rm {\bf
\Sigma} }  {\rm {\bf W}}^{ *} $, then ${\rm {\bf A}}^{ *}  = {\rm
{\bf W}}{\rm {\bf \Sigma} }^{ *}{\rm {\bf V}}^{ *} $ and ${\rm
{\bf A}}^{ +} = {\rm {\bf W}}{\rm {\bf \Sigma} }^{ +} {\rm {\bf
V}}^{ *} $. Since $ {\rm {\bf V}}$ is unitary,  then ${\rm {\bf
V}}^{ *}  = {\rm {\bf V}}^{ - 1}$. We have
\[\begin{array}{c}
  {\rm {\bf A}}^{ *} \left( {{\rm {\bf A}}{\rm {\bf A}}^{ *}  + \alpha {\rm
{\bf I}}} \right)^{ - 1} = {\rm {\bf W}}{\rm {\bf \Sigma} }{\rm
{\bf V}}^{ * } \cdot \left( {{\rm {\bf V}}{\rm {\bf \Sigma} }
\cdot {\rm {\bf W}}^{ * }{\rm {\bf W}} \cdot {\rm {\bf \Sigma}
}{\rm {\bf V}}^{ *}  + \alpha {\rm
{\bf I}}} \right)^{ - 1} =  \\
  ={\rm {\bf W}}{\rm {\bf \Sigma} }{\rm {\bf V}}^{ * } \cdot \left(
{{\rm {\bf V}}\left( {{\rm {\bf \Sigma} }{\rm {\bf \Sigma} }^{*} +
\alpha {\rm {\bf I}}} \right){\rm {\bf V}}^{ *} } \right)^{ - 1} =
{\rm {\bf W}}{\rm {\bf \Sigma} }\left( {{\rm {\bf \Sigma} }{\rm
{\bf \Sigma} }^{*} + \alpha {\rm {\bf I}}} \right)^{ - 1}{\rm {\bf
V}}^{ *}.
\end{array}
\]
 Consider the matrix
\[
{\rm {\bf \Sigma} }\left( {{\rm {\bf \Sigma} }{\rm {\bf \Sigma}
}^{*} + \alpha {\rm {\bf I}}} \right)^{ - 1} = \left(
{{\begin{array}{*{20}c}
 {{\frac{{\lambda _{1}} }{{\lambda _{1}^{2} + \alpha} }}} \hfill & {\ldots}
\hfill & {0} \hfill & {} \hfill & {} \hfill & {} \hfill \\
 {\ldots}  \hfill & {\ldots}  \hfill & {\ldots}  \hfill & {} \hfill & {{\rm
{\bf 0}}} \hfill & {} \hfill \\
 {0} \hfill & {\ldots}  \hfill & {{\frac{{\lambda _{r}} }{{\lambda _{r}^{2}
+ \alpha} }}} \hfill & {} \hfill & { \vdots}  \hfill & {} \hfill \\
 {} \hfill & { \vdots}  \hfill & {} \hfill & { \ddots}  \hfill & {} \hfill &
{} \hfill \\
 {} \hfill & {{\rm {\bf 0}}} \hfill & {} \hfill & {} \hfill & {{\rm {\bf
0}}} \hfill & {} \hfill \\
\end{array}} } \right).
\]
It is obviously that ${\mathop {\lim} \limits_{\alpha \to 0}} {\rm
{\bf \Sigma }}\left( {{\rm {\bf \Sigma} }{\rm {\bf \Sigma} }^{*} +
\alpha {\rm {\bf I}}} \right)^{ - 1} = {\rm {\bf \Sigma} }^{ +} $.
This implies ${\mathop {\lim} \limits_{\alpha \to 0} }{\rm {\bf
A}}^{ *} \left( {{\rm {\bf A}}{\rm {\bf A}}^{ *}  + \alpha {\rm
{\bf I}}} \right)^{ - 1}={\mathop {\lim} \limits_{\alpha \to 0}
}{\rm {\bf W}}{\rm {\bf \Sigma} }\left( {{\rm {\bf \Sigma} }{\rm
{\bf \Sigma} }^{*} + \alpha {\rm {\bf I}}} \right)^{ - 1}{\rm {\bf
V}}^{ *} = {\rm {\bf A}}^{ +} $.

By analogy we can prove that ${\rm {\bf A}}^{ +} ={\mathop {\lim}
\limits_{\alpha \to 0}} \left( {{\rm {\bf A}}^{
* }{\rm {\bf A}} + \alpha {\rm {\bf I}}} \right)^{ - 1}{\rm {\bf
A}}^{ *} $.$\blacksquare$

\begin{corollary}\label{kyrc14}  If
${\rm {\bf A}}\in {\mathbb H}^{m\times n} $, then the following
statements are true.
 \begin{itemize}
\item [ i)] If $\rm{rank}\,{\rm {\bf A}} = n$, then ${\rm {\bf A}}^{ +}
= \left( {{\rm {\bf A}}^{ *} {\rm {\bf A}}} \right)^{ - 1}{\rm
{\bf A}}^{ * }$ .
\item [ ii)] If $\rm{rank}\,{\rm {\bf A}} =
m$, then ${\rm {\bf A}}^{ +}  = {\rm {\bf A}}^{ * }\left( {{\rm
{\bf A}}{\rm {\bf A}}^{ *} } \right)^{ - 1}.$
\item [ iii)] If $\rm{rank}\,{\rm {\bf A}} = n = m$, then ${\rm {\bf
A}}^{ +}  = {\rm {\bf A}}^{ - 1}$ .
\end{itemize}
\end{corollary}

\section{ Determinantal representation of the Moore-Penrose
inverse.}
\newcommand{\rank}{\mathop{\rm rank}\nolimits}
\begin{lemma} \label{kyrc5}
If ${\rm {\bf A}} \in {\rm {\mathbb{H}}}^{m\times n}_r$, then $
 \rank\,\left( {{\rm {\bf A}}^{ *} {\rm {\bf A}}}
\right)_{.\,i} \left( {{\rm {\bf a}}_{.j}^{ *} }  \right) \le r. $
\end{lemma}
{\textit{Proof}}. Let's lead  elementary transformations of the
matrix $\left( {{\rm {\bf A}}^{ *} {\rm {\bf A}}} \right)_{.\,i}
\left( {{\rm {\bf a}}_{.j}^{ *} } \right)$ right-multiplying it by
elementary unimodular matrices ${\rm {\bf P}}_{i\,k} \left( { -
a_{jk}^{}} \right)$, $k \ne j$. The matrix ${\rm {\bf P}}_{\,i\,k}
\left( { - a_{jk}^{}} \right)$ has $-a_{j\,k} $ in the $(i, k)$
entry, 1 in all diagonal entries, and 0 in others. It is the
matrix of an elementary transformation. Right-multiplying  a
matrix by  ${\rm {\bf P}}_{\,i\,k} \left( { - a_{jk}^{}} \right)$,
where $k \ne j$, means adding to $k$-th column its $i$-th column
right-multiplying on $ - a_{jk} $. Then we get

\[
\left( {{\rm {\bf A}}^{ *} {\rm {\bf A}}} \right)_{.\,i} \left(
{{\rm {\bf a}}_{.\,j}^{ *} }  \right) \cdot {\prod\limits_{k \ne
i} {{\rm {\bf P}}_{i\,k} \left( {-a_{j\,k}}  \right) = {\mathop
{\left( {{\begin{array}{*{20}c}
 {{\sum\limits_{k \ne j} {a_{1k}^{ *}  a_{k1}} } } \hfill & {\ldots}  \hfill
& {a_{1j}^{ *} }  \hfill & {\ldots}  \hfill & {{\sum\limits_{k \ne
j } {a_{1k}^{ *}  a_{kn}}}} \hfill \\
 {\ldots}  \hfill & {\ldots}  \hfill & {\ldots}  \hfill & {\ldots}  \hfill &
{\ldots}  \hfill \\
 {{\sum\limits_{k \ne j} {a_{nk}^{ *}  a_{k1}} } } \hfill & {\ldots}  \hfill
& {a_{nj}^{ *} }  \hfill & {\ldots}  \hfill & {{\sum\limits_{k \ne
j } {a_{nk}^{ *}  a_{kn}}}} \hfill \\
\end{array}} }
\right)}\limits_{i-th}}}}.
\]

The obtained matrix  has the following factorization.

\[
{\mathop {\left( {{\begin{array}{*{20}c}
 {{\sum\limits_{k \ne j} {a_{1k}^{ *}  a_{k1}} } } \hfill & {\ldots}  \hfill
& {a_{1j}^{ *} }  \hfill & {\ldots}  \hfill & {{\sum\limits_{k \ne
j } {a_{1k}^{ *}  a_{kn}} } } \hfill \\
 {\ldots}  \hfill & {\ldots}  \hfill & {\ldots}  \hfill & {\ldots}  \hfill &
{\ldots}  \hfill \\
 {{\sum\limits_{k \ne j} {a_{nk}^{ *}  a_{k1}} } } \hfill & {\ldots}  \hfill
& {a_{nj}^{ *} }  \hfill & {\ldots}  \hfill & {{\sum\limits_{k \ne
j } {a_{nk}^{ *}  a_{kn}} } } \hfill \\
\end{array}} } \right)}\limits_{i-th}}  =
\]
\[
 = \left( {{\begin{array}{*{20}c}
 {a_{11}^{ *} }  \hfill & {a_{12}^{ *} }  \hfill & {\ldots}  \hfill &
{a_{1m}^{ *} }  \hfill \\
 {a_{21}^{ *} }  \hfill & {a_{22}^{ *} }  \hfill & {\ldots}  \hfill &
{a_{2m}^{ *} }  \hfill \\
 {\ldots}  \hfill & {\ldots}  \hfill & {\ldots}  \hfill & {\ldots}  \hfill
\\
 {a_{n1}^{ *} }  \hfill & {a_{n2}^{ *} }  \hfill & {\ldots}  \hfill &
{a_{nm}^{ *} }  \hfill \\
\end{array}} } \right){\mathop {\left( {{\begin{array}{*{20}c}
 {a_{11}}  \hfill & {\ldots}  \hfill & {0} \hfill & {\ldots}  \hfill &
{a_{n1}}  \hfill \\
 {\ldots}  \hfill & {\ldots}  \hfill & {\ldots}  \hfill & {\ldots}  \hfill &
{\ldots}  \hfill \\
 {0} \hfill & {\ldots}  \hfill & {1} \hfill & {\ldots}  \hfill & {0} \hfill
\\
 {\ldots}  \hfill & {\ldots}  \hfill & {\ldots}  \hfill & {\ldots}  \hfill &
{\ldots}  \hfill \\
 {a_{m1}}  \hfill & {\ldots}  \hfill & {0} \hfill & {\ldots}  \hfill &
{a_{mn}}  \hfill \\
\end{array}} } \right)}\limits_{i-th}} j-th.
\]
Denote by ${\rm {\bf \tilde {A}}}: = {\mathop {\left(
{{\begin{array}{*{20}c}
 {a_{11}}  \hfill & {\ldots}  \hfill & {0} \hfill & {\ldots}  \hfill &
{a_{n1}}  \hfill \\
 {\ldots}  \hfill & {\ldots}  \hfill & {\ldots}  \hfill & {\ldots}  \hfill &
{\ldots}  \hfill \\
 {0} \hfill & {\ldots}  \hfill & {1} \hfill & {\ldots}  \hfill & {0} \hfill
\\
 {\ldots}  \hfill & {\ldots}  \hfill & {\ldots}  \hfill & {\ldots}  \hfill &
{\ldots}  \hfill \\
 {a_{m1}}  \hfill & {\ldots}  \hfill & {0} \hfill & {\ldots}  \hfill &
{a_{mn}}  \hfill \\
\end{array}} } \right)}\limits_{i -th}} j -th$.
The matrix ${\rm {\bf \tilde {A}}}$ is obtained from ${\rm {\bf
A}}$ by replacing all entries of the $j$th row  and of the $i$th
column with zeroes except that the $(j, i)$ entry equals 1.
Elementary transformations of a matrix do not change its rank and
the rank of a matrix product does not exceed a rank of each
factor. It follows that $\rank\left( {{\rm {\bf A}}^{ *} {\rm {\bf
A}}} \right)_{.\,i} \left( {{\rm {\bf a}}_{.j}^{ *} } \right) \le
\min \,{\left\{ {\rank{\rm {\bf A}}^{ *}, \rank{\rm {\bf \tilde
{A}}}} \right\}}$. It is obvious that $\rank{\rm {\bf \tilde {A}}}
\ge \rank{\rm {\bf A}} = \rank{\rm {\bf A}}^{ *} $. Taking into
account Theorem \ref{kyrc6} we obtain $\rank{\rm {\bf A}}^{ *}
{\rm {\bf A}} = \rank{\rm {\bf A}}$. This completes the proof.
$\blacksquare$

The following lemma is  proved in the same way.
\begin{lemma}
If ${\rm {\bf A}} \in {\rm {\mathbb{H}}}^{m\times n}_r$, then
$\rank\left( {{\rm {\bf A}}{\rm {\bf A}}^{ *} } \right)_{.\,i}
\left( {{\rm {\bf a}}_{.j}^{ *} } \right) \le r $.
\end{lemma}
We shall use the following notations. Let $\alpha : = \left\{
{\alpha _{1} ,\ldots ,\alpha _{k}} \right\} \subseteq {\left\{
{1,\ldots ,m} \right\}}$ and $\beta : = \left\{ {\beta _{1}
,\ldots ,\beta _{k}} \right\} \subseteq {\left\{ {1,\ldots ,n}
\right\}}$ be subsets of the order $1 \le k \le \min {\left\{
{m,n} \right\}}$. By ${\rm {\bf A}}_{\beta} ^{\alpha} $ denote the
submatrix of ${\rm {\bf A}}$ determined by the rows indexed by
$\alpha$ and the columns indexed by $\beta$. Then ${\rm {\bf
A}}{\kern 1pt}_{\alpha} ^{\alpha}$ denotes the principal submatrix
determined by the rows and columns indexed by $\alpha$.
 If ${\rm {\bf A}} \in {\rm
M}\left( {n,{\rm {\mathbb{H}}}} \right)$ is Hermitian, then by
${\left| {{\rm {\bf A}}_{\alpha} ^{\alpha} } \right|}$ denote the
corresponding principal minor of $\det {\rm {\bf A}}$.
 For $1 \leq k\leq n$, denote by $\textsl{L}_{ k,
n}: = {\left\{ {\,\alpha :\alpha = \left( {\alpha _{1} ,\ldots
,\alpha _{k}} \right),\,{\kern 1pt} 1 \le \alpha _{1} \le \cdots
\le \alpha _{k} \le n} \right\}}$, the collection of strictly
increasing sequences of $k$ integers chosen from $\left\{
{1,\ldots ,n} \right\}$. For fixed $i \in \alpha $ and $j \in
\beta $, let $I_{r,\,m} {\left\{ {i} \right\}}: = {\left\{
{\,\alpha :\alpha \in L_{r,m} ,i \in \alpha}  \right\}}{\rm ,}
\quad J_{r,\,n} {\left\{ {j} \right\}}: = {\left\{ {\,\beta :\beta
\in L_{r,n} ,j \in \beta}  \right\}}$.

\begin{lemma}
If ${\rm {\bf A}} \in {\rm {\mathbb{H}}}^{m\times n}$ and $t \in
\mathbb{R}$, then

\begin{equation}
\label{kyr4} {\rm{cdet}} _{i} \left( {t{\rm {\bf I}} + {\rm {\bf
A}}^{ *} {\rm {\bf A}}} \right)_{.{\kern 1pt} i} \left( {{\rm {\bf
a}}_{.j}^{ *} }  \right) = c_{1}^{\left( {ij} \right)} t^{n - 1} +
c_{2}^{\left( {ij} \right)} t^{n - 2} + \cdots + c_{n}^{\left(
{ij} \right)},
\end{equation}

\noindent where $c_{n}^{\left( {ij} \right)} = {\rm{cdet}} _{i}
\left( {{\rm {\bf A}}^{ *} {\rm {\bf A}}} \right)_{.\,i} \left(
{{\rm {\bf a}}_{.\,j}^{ *} } \right)$  and $c_{k}^{\left( {ij}
\right)} = {\sum\limits_{\beta \in J_{k,\,n} {\left\{ {i}
\right\}}} {{\rm{cdet}} _{i} \left( {\left( {{\rm {\bf A}}^{ *}
{\rm {\bf A}}} \right)_{.\,i} \left( {{\rm {\bf a}}_{.\,j}^{ *} }
\right)} \right){\kern 1pt}  _{\beta} ^{\beta} } }$ for all $k
=1,\ldots,n - 1 $, $i = 1,\ldots,n$, and $j =1,\ldots,m$.
\end{lemma}
{\textit{Proof}}. Denote by ${\rm {\bf b}}_{.{\kern 1pt} {\kern
1pt} i} $ the $i$-th column of the Hermitian matrix ${\rm {\bf
A}}^{ *} {\rm {\bf A}} = :\left( {b_{ij}}\right)_{n\times n}$.
Consider the Hermitian matrix $\left( {t{\rm {\bf I}} + {\rm {\bf
A}}^{ *} {\rm {\bf A}}} \right)_{.{\kern 1pt} {\kern 1pt} i} ({\rm
{\bf b}}_{.{\kern 1pt} {\kern 1pt} i} ) \in {\rm
{\mathbb{H}}}^{n\times n}$. It differs from $\left( {t{\rm {\bf
I}} + {\rm {\bf A}}^{ *} {\rm {\bf A}}} \right)$ by an entry
$b_{ii} $. Taking into account Theorem \ref{kyrc7} we obtain
\begin{equation}
\label{kyr19} \det \left( {t{\rm {\bf I}} + {\rm {\bf A}}^{ *}
{\rm {\bf A}}} \right)_{.{\kern 1pt} i} \left( {{\rm {\bf
b}}_{.{\kern 1pt} {\kern 1pt} i}}  \right) = d_{1} t^{n - 1} +
d_{2} t^{n - 2} + \cdots + d_{n},
\end{equation}
where $d_{k} = {\sum\limits_{\beta \in J_{k,\,n} {\left\{ {i}
\right\}}} {\det \left( {{\rm {\bf A}}^{ *} {\rm {\bf A}}}
\right){\kern 1pt} {\kern 1pt} _{\beta} ^{\beta} } } $ is the sum
of all principal minors of order $k$ that contain the $i$-th
column for all $k = 1,\ldots,n - 1 $ and $d_{n} = \det \left(
{{\rm {\bf A}}^{ *} {\rm {\bf A}}} \right)$. Consequently we have
${\rm {\bf b}}_{.{\kern 1pt} {\kern 1pt} i} = \left(
{{\begin{array}{*{20}c}
 {{\sum\limits_{l} {a_{1l}^{ *}  a_{li}} } } \hfill \\
 {{\sum\limits_{l} {a_{2l}^{ *}  a_{li}} } } \hfill \\
 { \vdots}  \hfill \\
 {{\sum\limits_{l} {a_{nl}^{ *}  a_{li}} } } \hfill \\
\end{array}} } \right) = {\sum\limits_{l} {{\rm {\bf a}}_{.\,l}^{ *}  a_{li}}
}$, where ${\rm {\bf a}}_{.{\kern 1pt} {\kern 1pt} l}^{ *}  $ is
the $l$th column-vector  of ${\rm {\bf A}}^{ *}$ for all
$l=1,\ldots,m$. Taking into account Theorem \ref{kyrc2}, Lemma
\ref{kyrc1} and Proposition \ref{kyrc8} we obtain on the one hand
\begin{equation}
\label{kyr20}
\begin{array}{c}
  \det \left( {t{\rm {\bf I}} + {\rm {\bf A}}^{ *} {\rm {\bf A}}}
\right)_{.{\kern 1pt} i} \left( {{\rm {\bf b}}_{.{\kern 1pt}
{\kern 1pt} i} } \right) = {\rm{cdet}} _{i} \left( {t{\rm {\bf I}}
+ {\rm {\bf A}}^{ *} {\rm {\bf A}}} \right)_{.{\kern 1pt} i}
\left( {{\rm {\bf b}}_{.{\kern 1pt} {\kern 1pt} i}}  \right) =\\
   = {\sum\limits_{l} {{\rm{cdet}} _{i} \left( {t{\rm {\bf I}} + {\rm {\bf A}}^{ *
}{\rm {\bf A}}} \right)_{.{\kern 1pt} l} \left( {{\rm {\bf
a}}_{.{\kern 1pt} {\kern 1pt} l}^{ *}  a_{l{\kern 1pt} i}} \right)
= {\sum\limits_{l} {{\rm{cdet}} _{i} \left( {t{\rm {\bf I}} + {\rm
{\bf A}}^{ *} {\rm {\bf A}}} \right)_{.{\kern 1pt} i} \left( {{\rm
{\bf a}}_{.{\kern 1pt} {\kern 1pt} l}^{ *} }  \right) \cdot {\kern
1pt}} } }} a_{li}
\end{array}
\end{equation}
On the other hand having changed the order of summation,  we get
for all $k = 1,\ldots,n - 1 $
\begin{equation}
\label{kyr21}
\begin{array}{c}
  d_{k} = {\sum\limits_{\beta \in J_{k,\,n} {\left\{ {i} \right\}}}
{\det \left( {{\rm {\bf A}}^{ *} {\rm {\bf A}}} \right){\kern 1pt}
{\kern 1pt} _{\beta} ^{\beta} } }  = {\sum\limits_{\beta \in
J_{k,\,n} {\left\{ {i} \right\}}} {{\rm{cdet}} _{i} \left( {{\rm
{\bf A}}^{ *} {\rm {\bf A}}} \right){\kern 1pt} {\kern 1pt}
_{\beta} ^{\beta} } }  =\\
  {\sum\limits_{\beta \in J_{k,\,n}
{\left\{ {i} \right\}}} {{\sum\limits_{l} {{\rm{cdet}} _{i} \left(
{\left( {{\rm {\bf A}}^{ *} {\rm {\bf A}}{\kern 1pt}}
\right)_{.{\kern 1pt} {\kern 1pt} i} \left( {{\rm {\bf
a}}_{.{\kern 1pt} {\kern 1pt} l}^{ *} a_{l\,i}}  \right)}
\right)}} }} {\kern 1pt} _{\beta} ^{\beta}
 = {\sum\limits_{l} { {{\sum\limits_{\beta \in J_{k,\,n} {\left\{ {i}
\right\}}} {{\rm{cdet}} _{i} \left( {\left( {{\rm {\bf A}}^{ *}
{\rm {\bf A}}{\kern 1pt}}  \right)_{.{\kern 1pt} i} \left( {{\rm
{\bf a}}_{.{\kern 1pt} {\kern 1pt} l}^{ *} }  \right)}
\right){\kern 1pt} _{\beta} ^{\beta} } }} }}  \cdot a_{l{\kern
1pt} i}.
\end{array}
\end{equation}
By substituting (\ref{kyr20}) and (\ref{kyr21}) in (\ref{kyr19}),
and equating factors at $a _ {l \, i} $ when $l = j $, we obtain
the equality (\ref{kyr4}). $\blacksquare$

By analogy  the following lemma can be proved.
\begin{lemma}
If ${\rm {\bf A}} \in {\rm {\mathbb{H}}}^{m\times n}$ and $t \in
\mathbb{R}$, then
\[ {\rm{rdet}} _{j}  {( t{\rm {\bf I}} + {\rm
{\bf A}}{\rm {\bf A}}^{ *} )_{j\,.\,} ({\rm {\bf a}}_{i.}^{ *}  )}
   = r_{1}^{\left( {ij} \right)} t^{n
- 1} +r_{2}^{\left( {ij} \right)} t^{n - 2} + \cdots +
r_{n}^{\left( {ij} \right)},
\]

\noindent where  $r_{n}^{\left( {ij} \right)} = {\rm{rdet}} _{j}
{({\rm {\bf A}}{\rm {\bf A}}^{ *} )_{j\,.\,} ({\rm {\bf
a}}_{i.\,}^{ *} )}$ and $r_{k}^{\left( {ij} \right)} =
{{{\sum\limits_{\alpha \in I_{r,m} {\left\{ {j} \right\}}}
{{\rm{rdet}} _{j} \left( {({\rm {\bf A}}{\rm {\bf A}}^{ *}
)_{j\,.\,} ({\rm {\bf a}}_{i.\,}^{ *}  )} \right)\,_{\alpha}
^{\alpha} } }}}$ for all $k = 1,\ldots,n - 1 $, $i = 1,\ldots,n$,
and $j =1,\ldots,m$.
\end{lemma}
\begin{theorem}\label{kyrc16}
If ${\rm {\bf A}} \in {\rm {\mathbb{H}}}_{r}^{m\times n} $, then
the Moore-Penrose inverse  ${\rm {\bf A}}^{ +} = \left( {a_{ij}^{
+} } \right) \in {\rm {\mathbb{H}}}_{}^{n\times m} $ possess the
following determinantal representations:
\begin{equation}
\label{kyr5} a_{ij}^{ +}  = {\frac{{{\sum\limits_{\beta \in
J_{r,\,n} {\left\{ {i} \right\}}} {{\rm{cdet}} _{i} \left( {\left(
{{\rm {\bf A}}^{ *} {\rm {\bf A}}} \right)_{\,. \,i} \left( {{\rm
{\bf a}}_{.j}^{ *} }  \right)} \right){\kern 1pt} {\kern 1pt}
_{\beta} ^{\beta} } } }}{{{\sum\limits_{\beta \in J_{r,\,\,n}}
{{\left| {\left( {{\rm {\bf A}}^{ *} {\rm {\bf A}}} \right){\kern
1pt}  _{\beta} ^{\beta} }  \right|}}} }}},
\end{equation}
or
\begin{equation}
\label{kyr6} a_{ij}^{ +}  = {\frac{{{\sum\limits_{\alpha \in
I_{r,m} {\left\{ {j} \right\}}} {{\rm{rdet}} _{j} \left( {({\rm
{\bf A}}{\rm {\bf A}}^{ *} )_{j\,.\,} ({\rm {\bf a}}_{i.\,}^{ *}
)} \right)\,_{\alpha} ^{\alpha} } }}}{{{\sum\limits_{\alpha \in
I_{r,\,m}}  {{\left| {\left( {{\rm {\bf A}}{\rm {\bf A}}^{ *} }
\right){\kern 1pt}  _{\alpha} ^{\alpha} } \right|}}} }}}.
\end{equation}
\end{theorem}
{\textit{Proof}}.  At first we prove (\ref{kyr5}). Using Theorem
\ref{kyrc9}, we get ${\rm {\bf A}}^{ +}  = {\mathop {\lim}
\limits_{\alpha \to 0}} \left( {\alpha {\rm {\bf I}} + {\rm {\bf
A}}^{ *} {\rm {\bf A}}} \right)^{ - 1}{\rm {\bf A}}^{ * }$. The
matrix $\left( {\alpha {\rm {\bf I}} + {\rm {\bf A}}^{ *} {\rm
{\bf A}}} \right) \in {\rm {\mathbb{H}}}^{n\times n}$ is a
full-rank Hermitian matrix. Taking into account Theorem
\ref{kyrc10} it has an inverse, which we represent as a left
inverse matrix

\[
\left( {\alpha {\rm {\bf I}} + {\rm {\bf A}}^{ *} {\rm {\bf A}}}
\right)^{ - 1} = {\frac{{1}}{{\det \left( {\alpha {\rm {\bf I}} +
{\rm {\bf A}}^{ * }{\rm {\bf A}}} \right)}}}\left(
{{\begin{array}{*{20}c}
 {L_{11}}  \hfill & {L_{21}}  \hfill & {\ldots}  \hfill & {L_{n1}}  \hfill
\\
 {L_{12}}  \hfill & {L_{22}}  \hfill & {\ldots}  \hfill & {L_{n2}}  \hfill
\\
 {\ldots}  \hfill & {\ldots}  \hfill & {\ldots}  \hfill & {\ldots}  \hfill
\\
 {L_{1n}}  \hfill & {L_{2n}}  \hfill & {\ldots}  \hfill & {L_{nn}}  \hfill
\\
\end{array}} } \right),
\]
\noindent where $L_{ij} $ is a left $ij$th cofactor  of a matrix
$\alpha {\rm {\bf I}} + {\rm {\bf A}}^{ *} {\rm {\bf A}}$. Then we
have
\[\begin{array}{l}
   \left( {\alpha {\rm {\bf I}} + {\rm {\bf A}}^{ *} {\rm {\bf A}}} \right)^{ -
1}{\rm {\bf A}}^{ *}  = \\
  ={\frac{{1}}{{\det \left( {\alpha {\rm
{\bf I}} + {\rm {\bf A}}^{ *} {\rm {\bf A}}} \right)}}}\left(
{{\begin{array}{*{20}c}
 {{\sum\limits_{k = 1}^{n} {L_{k1} a_{k1}^{ *} } } } \hfill &
{{\sum\limits_{k = 1}^{n} {L_{k1} a_{k2}^{ *} } } } \hfill &
{\ldots}
\hfill & {{\sum\limits_{k = 1}^{n} {L_{k1} a_{km}^{ *} } } } \hfill \\
 {{\sum\limits_{k = 1}^{n} {L_{k2} a_{k1}^{ *} } } } \hfill &
{{\sum\limits_{k = 1}^{n} {L_{k2} a_{k2}^{ *} } } } \hfill &
{\ldots}
\hfill & {{\sum\limits_{k = 1}^{n} {L_{k2} a_{km}^{ *} } } } \hfill \\
 {\ldots}  \hfill & {\ldots}  \hfill & {\ldots}  \hfill & {\ldots}  \hfill
\\
 {{\sum\limits_{k = 1}^{n} {L_{kn} a_{k1}^{ *} } } } \hfill &
{{\sum\limits_{k = 1}^{n} {L_{kn} a_{k2}^{ *} } } } \hfill &
{\ldots}
\hfill & {{\sum\limits_{k = 1}^{n} {L_{kn} a_{km}^{ *} } } } \hfill \\
\end{array}} } \right).
\end{array}
\]
Using the definition of a left cofactor, we obtain
\begin{equation}
\label{kyr7} {\rm {\bf A}}^{ +}  = {\mathop {\lim} \limits_{\alpha
\to 0}} \left( {{\begin{array}{*{20}c}
 {{\frac{{{\rm cdet} _{1} \left( {\alpha {\rm {\bf I}} + {\rm {\bf A}}^{ *} {\rm
{\bf A}}} \right)_{.1} \left( {{\rm {\bf a}}_{.1}^{ *} }
\right)}}{{\det \left( {\alpha {\rm {\bf I}} + {\rm {\bf A}}^{ *}
{\rm {\bf A}}} \right)}}}} \hfill & {\ldots}  \hfill &
{{\frac{{{\rm cdet} _{1} \left( {\alpha {\rm {\bf I}} + {\rm {\bf
A}}^{ *} {\rm {\bf A}}} \right)_{.1} \left( {{\rm {\bf a}}_{.m}^{
*} }  \right)}}{{\det \left( {\alpha {\rm {\bf I}} + {\rm {\bf
A}}^{ *} {\rm {\bf A}}} \right)}}}} \hfill \\
 {\ldots}  \hfill & {\ldots}  \hfill & {\ldots}  \hfill \\
 {{\frac{{{\rm cdet} _{n} \left( {\alpha {\rm {\bf I}} + {\rm {\bf A}}^{ *} {\rm
{\bf A}}} \right)_{.n} \left( {{\rm {\bf a}}_{.1}^{ *} }
\right)}}{{\det \left( {\alpha {\rm {\bf I}} + {\rm {\bf A}}^{ *}
{\rm {\bf A}}} \right)}}}} \hfill & {\ldots}  \hfill &
{{\frac{{{\rm cdet} _{n} \left( {\alpha {\rm {\bf I}} + {\rm {\bf
A}}^{ *} {\rm {\bf A}}} \right)_{.n} \left( {{\rm {\bf a}}_{.m}^{
*} }  \right)}}{{\det \left( {\alpha {\rm {\bf I}} + {\rm {\bf
A}}^{ *} {\rm {\bf A}}} \right)}}}} \hfill \\
\end{array}} } \right).
\end{equation}
By Theorem \ref{kyrc7} we have $\det \left( {\alpha {\rm {\bf I}}
+ {\rm {\bf A}}^{ *} {\rm {\bf A}}} \right) = \alpha ^{n} + d_{1}
\alpha ^{n - 1} + d_{2} \alpha ^{n - 2} + \cdots + d_{n} $, where
$d_{k} = {\sum\limits_{\beta \in J_{k,\,n}}  {{\left| {\left(
{{\rm {\bf A}}^{ *} {\rm {\bf A}}} \right){\kern 1pt} {\kern 1pt}
_{\beta} ^{\beta} }  \right|}}} $ is a sum of principal minors of
${\rm {\bf A}}^{ *} {\rm {\bf A}}$ of order $k$ for all  $k =
1,\ldots,n - 1 $ and $d_{n} = \det {\rm {\bf A}}^{ *} {\rm {\bf
A}}$. Since $\rank{\rm {\bf A}}^{ *} {\rm {\bf A}} = \rank{\rm
{\bf A}} = r$ and $d_{n} = d_{n - 1} = \cdots = d_{r + 1} = 0$, it
follows that $\det \left( {\alpha {\rm {\bf I}} + {\rm {\bf A}}^{
*} {\rm {\bf A}}} \right) = \alpha ^{n} + d_{1} \alpha ^{n - 1} +
d_{2} \alpha ^{n - 2} + \cdots + d_{r} \alpha ^{n - r}$. Using
(\ref{kyr4})  we get ${\rm{cdet}} _{i} \left( {\alpha {\rm {\bf
I}} + {\rm {\bf A}}^{ * }{\rm {\bf A}}} \right)_{.i} \left( {{\rm
{\bf a}}_{.j}^{ *} }  \right) = c_{1}^{\left( {ij} \right)} \alpha
^{n - 1} + c_{2}^{\left( {ij} \right)} \alpha ^{n - 2} + \cdots +
c_{n}^{\left( {ij} \right)} $ for all $i =1,\ldots,n $ and $j
=1,\ldots,m $, where $c_{k}^{\left( {ij} \right)} =
{\sum\limits_{\beta \in J_{k,\,n} {\left\{ {i} \right\}}}
{{\rm{cdet}} _{i} \left( {({\rm {\bf A}}^{ *} {\rm {\bf
A}})_{.\,i} \left( {{\rm {\bf a}}_{.j}^{ *} } \right)}
\right){\kern 1pt} {\kern 1pt} _{\beta }^{\beta} } } $ for all $k
= 1,\ldots,n - 1 $ and $c_{n}^{\left( {ij} \right)} = {\rm{cdet}}
_{i} \left( {{\rm {\bf A}}^{ *} {\rm {\bf A}}} \right)_{.i} \left(
{{\rm {\bf a}}_{.j}^{ *} }  \right)$.

Now we prove that $c_{k}^{\left( {ij} \right)} = 0$, when $k \ge r
+ 1$ for all $i = 1,\ldots,n $, and $j =1,\ldots,m $. By Lemma
\ref{kyrc5} $\rank\left( {{\rm {\bf A}}^{ * }{\rm {\bf A}}}
\right)_{.\,i} \left( {{\rm {\bf a}}_{.j}^{ *} } \right) \le r$,
then the matrix $\left( {{\rm {\bf A}}^{ *} {\rm {\bf A}}}
\right)_{.\,i} \left( {{\rm {\bf a}}_{.j}^{ *} } \right)$ has no
more $r$ right-linearly independent columns.

Consider $\left( {({\rm {\bf A}}^{ * }{\rm {\bf A}})_{\,.\,i}
\left( {{\rm {\bf a}}_{.j}^{ *} } \right)} \right){\kern 1pt}
{\kern 1pt} _{\beta} ^{\beta}  $, when $\beta \in J_{k,n} {\left\{
{i} \right\}}$. It is a principal submatrix of  $\left( {{\rm {\bf
A}}^{ *} {\rm {\bf A}}} \right)_{.\,i} \left( {{\rm {\bf
a}}_{.j}^{ *} } \right)$ of order $k \ge r + 1$. Deleting both its
$i$-th row and column, we obtain a principal submatrix of order $k
- 1$ of  ${\rm {\bf A}}^{ *} {\rm {\bf A}}$.
 We denote it by ${\rm {\bf M}}$. The following cases are
possible.

Let $k = r + 1$ and $\det {\rm {\bf M}} \ne 0$. In this case all
columns of $ {\rm {\bf M}} $ are right-linearly independent. The
addition of all of them on one coordinate to columns of
 $\left( {\left( {{\rm {\bf A}}^{ *} {\rm {\bf A}}}
\right)_{.\,i} \left( {{\rm {\bf a}}_{.j}^{ *} }  \right)} \right)
{\kern 1pt} _{\beta} ^{\beta}$ keeps their right-linear
independence. Hence, they are basis in a matrix $\left( {\left(
{{\rm {\bf A}}^{ *} {\rm {\bf A}}} \right)_{\,.\,i} \left( {{\rm
{\bf a}}_{.j}^{ *} } \right)} \right){\kern 1pt} {\kern 1pt}
_{\beta} ^{\beta} $, and by Theorem \ref{kyrc11} the $i$-th column
is the right linear combination of its basis columns. From this by
Theorem \ref{kyrc12}, we get ${\rm{cdet}} _{i} \left( {\left(
{{\rm {\bf A}}^{ *} {\rm {\bf A}}} \right)_{\,.\,i} \left( {{\rm
{\bf a}}_{.j}^{ *} } \right)} \right){\kern 1pt} {\kern 1pt}
_{\beta} ^{\beta}  = 0$, when $\beta \in J_{k,n} {\left\{ {i}
\right\}}$ and $k \ge r + 1$.

If $k = r + 1$ and $\det {\rm {\bf M}} = 0$, then $p$, ($p < k$),
columns are basis in  ${\rm {\bf M}}$ and in  $\left( {\left(
{{\rm {\bf A}}^{ *} {\rm {\bf A}}} \right)_{.\,i} \left( {{\rm
{\bf a}}_{.j}^{ *} } \right)} \right){\kern 1pt} {\kern 1pt}
_{\beta} ^{\beta} $. Then by Theorems \ref{kyrc11} and
\ref{kyrc12} we obtain ${\rm{cdet}} _{i} \left( {\left( {{\rm {\bf
A}}^{ *} {\rm {\bf A}}} \right)_{\,.\,i} \left( {{\rm {\bf
a}}_{.j}^{ *} } \right)} \right){\kern 1pt} {\kern 1pt} _{\beta}
^{\beta}  = 0$ as well.

If $k > r + 1$, then from Theorems \ref{kyrc6} and \ref{kyrc13} it
follows that $\det {\rm {\bf M}} = 0$ and $p$, ($p < k - 1$),
columns are basis in the both matrices ${\rm {\bf M}}$ and $\left(
{\left( {{\rm {\bf A}}^{ *} {\rm {\bf A}}} \right)_{\,.\,i} \left(
{{\rm {\bf a}}_{.j}^{ *} }  \right)} \right){\kern 1pt} {\kern
1pt} _{\beta} ^{\beta}  $. Then by Theorems \ref{kyrc11} and
\ref{kyrc12}, we obtain that ${\rm{cdet}} _{i} \left( {\left(
{{\rm {\bf A}}^{ *} {\rm {\bf A}}} \right)_{\,.\,i} \left( {{\rm
{\bf a}}_{.j}^{ *} } \right)} \right){\kern 1pt} {\kern 1pt}
_{\beta} ^{\beta}  = 0$.

Thus in all cases we have ${\rm{cdet}} _{i} \left( {\left( {{\rm
{\bf A}}^{ * }{\rm {\bf A}}} \right)_{\,.\,i} \left( {{\rm {\bf
a}}_{.j}^{ *} }  \right)} \right){\kern 1pt} {\kern 1pt} _{\beta}
^{\beta}  = 0$, when $\beta \in J_{k,n} {\left\{ {i} \right\}}$
and $r + 1 \le k < n$. From here if $r + 1 \le k < n$, then
$c_{k}^{\left( {ij} \right)} = {\sum\limits_{\beta \in J_{k,\,n}
{\left\{ {i} \right\}}} {{\rm{cdet}} _{i} \left( {\left( {{\rm
{\bf A}}^{ *} {\rm {\bf A}}} \right)_{\,.\,i} \left( {{\rm {\bf
a}}_{.j}^{ *} } \right)} \right) {\kern 1pt} _{\beta} ^{\beta} }
}=0 $, and $c_{n}^{\left( {ij} \right)} = {\rm{cdet}} _{i} \left(
{{\rm {\bf A}}^{ *} {\rm {\bf A}}} \right)_{.\,i} \left( {{\rm
{\bf a}}_{.j}^{ *} }  \right) = 0$ for all $i = 1,\ldots,n $ and
$j = 1,\ldots,m $.

Hence, ${\rm{cdet}} _{i} \left( {\alpha {\rm {\bf I}} + {\rm {\bf
A}}^{ *} {\rm {\bf A}}} \right)_{.\,i} \left( {{\rm {\bf
a}}_{.\,j}^{ *} }  \right) =c_{1}^{\left( {ij} \right)} \alpha ^{n
- 1} + c_{2}^{\left( {ij} \right)} \alpha ^{n - 2} + \cdots +
c_{r}^{\left( {ij} \right)} \alpha ^{n - r}$ for all $i =
1,\ldots,n $ and $j = 1,\ldots,m $. By substituting these values
in the matrix from (\ref{kyr7}), we obtain
\[\begin{array}{c}
  {\rm {\bf A}}^{ +}  = {\mathop {\lim} \limits_{\alpha \to 0}} \left(
{{\begin{array}{*{20}c}
 {{\frac{{c_{1}^{\left( {11} \right)} \alpha ^{n - 1} + \cdots +
c_{r}^{\left( {11} \right)} \alpha ^{n - r}}}{{\alpha ^{n} + d_{1}
\alpha ^{n - 1} + \cdots + d_{r} \alpha ^{n - r}}}}} \hfill &
{\ldots}  \hfill & {{\frac{{c_{1}^{\left( {1m} \right)} \alpha ^{n
- 1} + \cdots + c_{r}^{\left( {1m} \right)} \alpha ^{n -
r}}}{{\alpha ^{n} + d_{1} \alpha
^{n - 1} + \cdots + d_{r} \alpha ^{n - r}}}}} \hfill \\
 {\ldots}  \hfill & {\ldots}  \hfill & {\ldots}  \hfill \\
 {{\frac{{c_{1}^{\left( {n1} \right)} \alpha ^{n - 1} + \cdots +
c_{r}^{\left( {n1} \right)} \alpha ^{n - r}}}{{\alpha ^{n} + d_{1}
\alpha ^{n - 1} + \cdots + d_{r} \alpha ^{n - r}}}}} \hfill &
{\ldots}  \hfill & {{\frac{{c_{1}^{\left( {nm} \right)} \alpha ^{n
- 1} + \cdots + c_{r}^{\left( {nm} \right)} \alpha ^{n -
r}}}{{\alpha ^{n} + d_{1} \alpha
^{n - 1} + \cdots + d_{r} \alpha ^{n - r}}}}} \hfill \\
\end{array}} } \right) =\\
  \left( {{\begin{array}{*{20}c}
 {{\frac{{c_{r}^{\left( {11} \right)}} }{{d_{r}} }}} \hfill & {\ldots}
\hfill & {{\frac{{c_{r}^{\left( {1m} \right)}} }{{d_{r}} }}} \hfill \\
 {\ldots}  \hfill & {\ldots}  \hfill & {\ldots}  \hfill \\
 {{\frac{{c_{r}^{\left( {n1} \right)}} }{{d_{r}} }}} \hfill & {\ldots}
\hfill & {{\frac{{c_{r}^{\left( {nm} \right)}} }{{d_{r}} }}} \hfill \\
\end{array}} } \right).
\end{array}
\]

Here $c_{r}^{\left( {ij} \right)} = {\sum\limits_{\beta \in
J_{r,\,n} {\left\{ {i} \right\}}} {{\rm{cdet}} _{i} \left( {\left(
{{\rm {\bf A}}^{ *} {\rm {\bf A}}} \right)_{\,.\,i} \left( {{\rm
{\bf a}}_{.j}^{ *} }  \right)} \right){\kern 1pt} {\kern 1pt}
_{\beta} ^{\beta} } } $ and $d_{r} = {\sum\limits_{\beta \in
J_{r,\,\,n}} {{\left| {\left( {{\rm {\bf A}}^{ * }{\rm {\bf A}}}
\right){\kern 1pt} {\kern 1pt} _{\beta} ^{\beta} } \right|}}} $.
Thus, we have obtained the determinantal representation of ${\rm
{\bf A}}_{}^{ +}  $  (\ref{kyr5}).

Similarly one can prove the determinantal representation of ${\rm
{\bf A}}_{}^{ +} $  (\ref{kyr6}). $\blacksquare$
\begin{remark} In (\ref{kyr5}) the index $i$ in ${{\rm{cdet}} _{i} \left( {\left(
{{\rm {\bf A}}^{ *} {\rm {\bf A}}} \right)_{\,.{\kern 1pt} \,i}
\left( {{\rm {\bf a}}_{.j}^{ *} }  \right)} \right){\kern 1pt}
{\kern 1pt} _{\beta} ^{\beta} }$ designates $i$th column of
$\left( {\left( {{\rm {\bf A}}^{ *} {\rm {\bf A}}}
\right)_{\,.{\kern 1pt} \,i} \left( {{\rm {\bf a}}_{.j}^{ *} }
\right)} \right)$, but in the submatrix $\left( {\left( {{\rm {\bf
A}}^{ *} {\rm {\bf A}}} \right)_{\,.{\kern 1pt} \,i} \left( {{\rm
{\bf a}}_{.j}^{ *} }  \right)} \right){\kern 1pt} {\kern 1pt}
_{\beta} ^{\beta} $ the entries of ${{\rm {\bf a}}_{.j}^{ *} }$
may be placed in a column with the another index. In (\ref{kyr6})
we have equivalently.
\end{remark}

\begin{remark}\label{kyrc19}
If $\rank{\rm {\bf A}} = n$, then by Corollary \ref{kyrc14} ${\rm
{\bf A}}^{ +}  = \left( {{\rm {\bf A}}^{ *} {\rm {\bf A}}}
\right)^{ - 1}{\rm {\bf A}}^{ *} $. Considering $\left( {{\rm {\bf
A}}^{ *} {\rm {\bf A}}} \right)^{ - 1}$ as a left inverse, we get
the following representation of  ${\rm {\bf A}}^{ +} $:

\begin{equation}
\label{kyr8}{\rm {\bf A}}^{ +}  = {\frac{{1}}{{\rm {ddet} {\rm
{\bf A}}}}}
\begin{pmatrix}
  {{\rm{cdet}} _{1} ({{\rm {\bf A}}^{ *} {\rm {\bf
A}}})_{.\,1} \left( {{\rm {\bf a}}_{.\,1}^{ *} }  \right)} &
\ldots & {{\rm{cdet}} _{1} ({{\rm {\bf A}}^{ *} {\rm {\bf
A}}})_{.\,1} \left(
{{\rm {\bf a}}_{.\,m}^{ *} }  \right)} \\
  \ldots & \ldots & \ldots \\
  {{\rm{cdet}} _{n} ({{\rm {\bf A}}^{ *} {\rm {\bf A}}})_{.\,n}
 \left( {{\rm {\bf a}}_{.\,1}^{ *} }  \right)} & \ldots &
 {{\rm{cdet}} _{n} ({{\rm {\bf A}}^{ *} {\rm {\bf
A}}})_{.\,n} \left( {{\rm {\bf a}}_{.\,m}^{ *} }  \right)}.
\end{pmatrix}.
\end{equation}
If $m > n$, then by Theorem \ref{kyrc16} for ${\rm {\bf A}}^{ +} $
we have (\ref{kyr5}) as well.
\end{remark}
\begin{remark}\label{kyrc17}
If $\rank{\rm {\bf A}} = m$, then by Corollary \ref{kyrc14} ${\rm
{\bf A}}^{ +}  = {\rm {\bf A}}^{ *} \left( {{\rm {\bf A}}{\rm {\bf
A}}^{ *} } \right)^{ - 1}$. Considering  $\left( {{\rm {\bf
A}}{\rm {\bf A}}^{ *} } \right)^{ - 1}$ as a right inverse, we get
the following representation of ${\rm {\bf A}}^{ +} $:

\begin{equation}
\label{kyr9} {\rm {\bf A}}^{ +}  = {\frac{{1}}{{\rm {ddet} {{\rm
{\bf A}} }}}}
\begin{pmatrix}
 {{\rm{rdet}} _{1} ({\rm {\bf A}} {\rm {\bf A}}^{*})_{1.} \left( {{\rm {\bf a}}_{1.}^{ *} }  \right)} & \ldots &
 {{\rm{rdet}} _{m} ({\rm {\bf A}} {\rm {\bf A}}^{*})_{m.} \left( {{\rm {\bf a}}_{1.}^{ *} }  \right)}\\
 \ldots & \ldots & \ldots \\
  {{\rm{rdet}} _{1} ({\rm {\bf A}} {\rm {\bf A}}^{*})_{1.} \left( {{\rm {\bf a}}_{n.}^{ *} }  \right)} & \ldots &
  {{\rm{rdet}} _{m} ({\rm {\bf A}} {\rm {\bf A}}^{*})_{m\,.}
\left( {{\rm {\bf a}}_{n\,.}^{ *} }  \right)}
\end{pmatrix}.
\end{equation}
If $m < n$, then by Theorem \ref{kyrc16} for ${\rm {\bf A}}^{ +} $
we also have (\ref{kyr6}).

\end{remark}

\begin{corollary}
If ${\rm {\bf A}} \in {\rm {\mathbb{H}}}_{r}^{m\times n} $, where
$r < \min {\left\{ {m,n} \right\}}$ or $r = m < n$, then for a
projection matrix ${\rm {\bf A}}^{ +} {\rm {\bf A}} = :{\rm {\bf
P}} = \left( {p_{ij}} \right)_{n\times n} $ we  have its following
determinantal representation
\[
p_{ij} = {\frac{{{\sum\limits_{\beta \in J_{r,\,\,n} {\left\{ {i}
\right\}}} {{\rm{cdet}} _{i} \left( {\left( {{\rm {\bf A}}^{ *}
{\rm {\bf A}}} \right)_{.\,i} \left({\rm {\bf d}}_{.j} \right)}
\right){\kern 1pt}  _{\beta} ^{\beta} } }}}{{{\sum\limits_{\beta
\in J_{r,\,n}}  {{\left| {\left( {{\rm {\bf A}}^{ *} {\rm {\bf
A}}} \right){\kern 1pt}  _{\beta} ^{\beta} } \right|}}} }}},
\]

\noindent where ${\rm {\bf d}}_{.j} $ is the $j$-th column of
${{\rm {\bf A}}^{ *} {\rm {\bf A}}} \in {\rm
{\mathbb{H}}}^{n\times n}$ and for all $i,j = 1,\ldots,n$.

\end{corollary}
{\textit{Proof}}. Representing ${\rm {\bf A}}^{ +} $ by (\ref
{kyr5}) and right-multiplying  it
  by $ {\rm {\bf A}} $, we obtain for  an entry $p _ {ij} $ of $ {\rm {\bf
A}}^{+} {\rm {\bf A}} =: {\rm {\bf P}} = \left ({p _ {ij}} \right)
_ {n \times n} $.

\[\begin{array}{c}
  p_{ij} = {\sum\limits_{p = 1}^{m} {c_{ip} \cdot a_{pj}} }  =
{\sum\limits_{k} {{\frac{{{\sum\limits_{\beta \in J_{r,\,n}
{\left\{ {i} \right\}}} {{\rm{cdet}} _{i} \left( {\left( {{\rm
{\bf A}}^{ *} {\rm {\bf A}}} \right)_{.\,i} \left({\rm {\bf
a}}_{.\,j}^{ *}\right) }  \right) {\kern 1pt} _{\beta} ^{\beta} }
} }}{{{\sum\limits_{\beta \in J_{r,\,n}}  {{\left| {\left( {{\rm
{\bf A}}^{ *} {\rm {\bf A}}} \right) {\kern 1pt}
_{\beta} ^{\beta} }  \right|}}} }}}}}  \cdot a_{kj} =\\
   = {\frac{{{\sum\limits_{\beta \in J_{r,\,n} {\left\{ {i} \right\}}}
{{\sum\limits_{k} {{\rm{cdet}} _{i} \left( {\left( {{\rm {\bf
A}}^{ *} {\rm {\bf A}}} \right)_{\,.\,i} \left({\rm {\bf
a}}_{.j}^{ *}\right) } \right){\kern 1pt} {\kern 1pt} _{\beta}
^{\beta} } }  \cdot \,a_{kj}} } }}{{{\sum\limits_{\beta \in
J_{r,\,\,n}}  {{\left| {\left( {{\rm {\bf A}}^{ *} {\rm {\bf A}}}
\right) {\kern 1pt} _{\beta} ^{\beta} }  \right|}}} }}} =
{\frac{{{\sum\limits_{\beta \in J_{r,\,n} {\left\{ {i} \right\}}}
{{\rm{cdet}} _{i} \left( {\left( {{\rm {\bf A}}^{ *} {\rm {\bf
A}}} \right)_{.\,i} \left({ {\rm {\bf d}}}_{.\,j}\right)}  \right)
{\kern 1pt} _{\beta} ^{\beta} } }}}{{{\sum\limits_{\beta \in
J_{r,\,n}} {{\left| {\left( {{\rm {\bf A}}^{ *} {\rm {\bf A}}}
\right) {\kern 1pt} _{\beta} ^{\beta} } \right|}}} }}},
\end{array}
\]
\noindent where $d_{.j} $ is the $j$th column of ${{\rm {\bf A}}^{
*} {\rm {\bf A}}} \in {\rm {\mathbb{H}}}_{}^{n\times n}$ and for
all $i, j =1,\ldots,n $. $\blacksquare$

By analogy can be proved the following corollary.
\begin{corollary}
If ${\rm {\bf A}} \in {\rm {\mathbb{H}}}_{r}^{m\times n}$, where
$r < \min {\left\{ {m,n} \right\}}$ or $r = n < m$, then for the
projection matrix ${\rm {\bf A}}{\rm {\bf A}}^{ +} = :{\rm {\bf
Q}} = \left( {q_{ij}} \right)_{m\times m} $ we  have its following
determinantal representation
\[
q_{ij} = {\frac{{{\sum\limits_{\alpha \in I_{r,\,\,m} {\left\{ {i}
\right\}}} {{\left| {\left( {({\rm {\bf A}}\,{\rm {\bf A}}^{ *}
)_{i{\kern 1pt} .}\, ({\rm {\bf g}}  _{j{\kern 1pt}  .}\, )}
\right) {\kern 1pt} _{\alpha} ^{\alpha} } \right|}}}
}}{{{\sum\limits_{\alpha \in I_{r,\,m}} {{\left| {\left( {{\rm
{\bf A}}{\rm {\bf A}}^{ *} } \right){\kern 1pt} _{\alpha
}^{\alpha} }  \right|}}} }}},
\]
\noindent where ${\rm {\bf g}}_{j.} $ is the $j$th row of $({\rm
{\bf A}} {\rm {\bf A}}^{*})\in {\rm {\mathbb{H}}}^{m\times m}$ and
for all $i,j = 1,\ldots,m$.

\end{corollary}
\begin{remark}
By definition of a classical adjoint matrix of ${\rm {\bf A}} \in
{\rm {\mathbb{C}}}^{n\times n}$,  ${\rm Adj}\,{\left[ {{\rm {\bf
A}}} \right]} \cdot {\rm {\bf A}} = {\rm {\bf A}} \cdot {\rm
Adj}\,{\left[ {{\rm {\bf A}}} \right]} = \det {\rm {\bf A}} \cdot
{\rm {\bf I}}$. Let ${\rm {\bf A}} \in \mathbb{H}^{m\times n}$. If
$\rank {\rm {\bf A}} = n$, the by Corollary \ref{kyrc14} we have
${\rm {\bf A}}^{ +} {\rm {\bf A}} = {\rm {\bf I}}_{n}$.
Representing ${\rm {\bf A}}^{ +}$ by (\ref{kyr8}) as ${\rm {\bf
A}}^{ +}  = {\frac{{{\rm {\bf L}}}}{{\det \left( {{\rm {\bf A}}^{
*} {\rm {\bf A}}} \right)}}}$, where ${\rm {\bf L}} = \left(
{{\rm{cdet}} _{i} \left( {\left( {{\rm {\bf A}}^{ *} {\rm {\bf
A}}} \right)_{.\,i} \left( {\rm {\bf a}}_{.j}^{ *} \right)  {\kern
1pt}} \right)} \right)_{n\times m}$, we obtain ${\rm {\bf L}}{\rm
{\bf A}} = \det \left( {{\rm {\bf A}}^{ *} {\rm {\bf A}}} \right)
\cdot {\rm {\bf I}}_{n} $. This means that the matrix ${\rm {\bf
L}}=: {\rm Adj}\, _{L}  {\kern 1pt} {\left[ {{\rm {\bf A}}}
\right]} $ is the left classical adjoint matrix of ${\rm {\bf A}}
\in \mathbb{H}^{m\times n}$.

If $\rank {\rm {\bf A}} = m$, then by definition of a right
classical adjoint matrix of ${\rm {\bf A}} \in {\rm
{\mathbb{H}}}^{m\times n}$ by Corollary \ref{kyrc14} and by
(\ref{kyr9}) we can put

\[{\rm Adj}\, _{R}  {\kern 1pt} {\left[ {{\rm {\bf A}}} \right]}: =
\left( {\left( {{\rm{rdet}} _{j} ({\rm {\bf A}}{\rm {\bf A}}^{ *}
)_{j.} {\left({\rm {\bf a}}_{i.}^{ *}\right) } }\right)\,}
\right)_{m\times n}, \]

\noindent since in this case ${\rm {\bf A}} \cdot {\rm Adj}\, _{R}
{\kern 1pt} {\left[ {{\rm {\bf A}}} \right]} = \det ({\rm {\bf
A}}{\rm {\bf A}}^{ *} ) \cdot {\rm {\bf I}}$.

If $\rank {\rm {\bf A}} = r < \min {\left\{ {m,n} \right\}}$, then
an analog of a left classical adjoint matrix of ${\rm {\bf A}} \in
{\rm {\mathbb{H}}}^{m\times n}$ by (\ref{kyr5}) can accept
\[{\rm
Adj}\, _{L} {\kern 1pt} {\kern 1pt} {\left[ {{\rm {\bf A}}}
\right]}: = \left( {{\kern 1pt} {\sum\limits_{\alpha \in J_{r,\,n}
{\left\{ {i} \right\}}} {{\rm{cdet}} _{i} \left( {({\rm {\bf A}}^{
* }{\rm {\bf A}})_{.\,i} \,{\rm {\bf a}}_{.\,j}^{ *} }
\right){\kern 1pt}  _{\alpha} ^{\alpha} } } {\kern 1pt}}
\right)_{n\times m}. \]
\noindent Indeed, since eigenvalues of a
projection matrix are only 1 and 0, there exists such a unitary
matrix ${\rm {\bf U}}\in {\rm {\mathbb{H}}}^{n\times n}$ that
\[\begin{array}{c}
  {\rm Adj}\,{\kern 1pt} _{L}  {\kern 1pt} {\left[ {{\rm {\bf A}}} \right]}
\cdot {\rm {\bf A}} = {\sum\limits_{\alpha \in I_{r,\,n}} {{\left|
{\left( {{\rm {\bf A}}^{ *} {\rm {\bf A}}} \right) {\kern 1pt}
_{\alpha} ^{\alpha} }  \right|}}}  \cdot {\rm {\bf P}}=\\
   = {\sum\limits_{\alpha \in I_{r,\,n}}  {{\left| {\left( {{\rm {\bf A}}^{
*} {\rm {\bf A}}} \right) {\kern 1pt} _{\alpha} ^{\alpha} }
\right|}}}  \cdot {\rm {\bf U}}{\rm diag}(1,\ldots ,1,0,\ldots
,0){\rm {\bf U}}^{ *}. \end{array}
\]

If $\rank {\rm {\bf A}} = r < \min {\left\{ {m,n} \right\}}$, then
by an analogue  of a right classical adjoint matrix of ${\rm {\bf
A}} \in {\rm {\mathbb{H}}}^{m\times n}$ by (\ref{kyr6})  we can
put
\[{\rm Adj}\,_{R} {\kern 1pt} {\kern 1pt} {\left[ {{\rm {\bf
A}}} \right]}: = \left( {\,{\sum\limits_{\alpha \in I_{r,\,m}
{\left\{ {j} \right\}}} {{\rm{rdet}} _{j} \left( {({\rm {\bf
A}}{\rm {\bf A}}^{ *} )_{j.} \left({\rm {\bf a}}_{i.}^{ *}\right)
{\kern 1pt}} \right)_{\alpha} ^{\alpha} } } {\kern 1pt}}
\right)_{n\times m},\]
\noindent as there exists such a unitary
matrix ${\rm {\bf V}}\in {\rm {\mathbb{H}}}^{m\times m}$ that
\[\begin{array}{c}
  {\rm {\bf A}} \cdot {\rm Adj}\, _{R} {\kern 1pt} {\kern 1pt}
{\left[ {{\rm {\bf A}}} \right]} = {\sum\limits_{\alpha \in
J_{r,\,m}}  {{\left| {\left( {{\rm {\bf A}} {\rm {\bf A}}^{ *}}
\right) {\kern 1pt} _{\alpha} ^{\alpha} }  \right|}}} \cdot
{\rm {\bf Q}}=\\
 = {\sum\limits_{\alpha \in J_{r,\,m}}  {{\left|
{\left( {{\rm {\bf A}} {\rm {\bf A}}^{ *}} \right) {\kern 1pt}
_{\alpha} ^{\alpha} }  \right|}}}  \cdot {\rm {\bf V}}{\rm
diag}(1,\ldots ,1,0,\ldots ,0){\rm {\bf V}}^{ *}.
\end{array}
\]
\end{remark}

\begin{remark} If ${\rm {\bf A}} \in {\mathbb C}^{m\times n}$ is a matrix with
complex entries, then we obtain the following analogs of
(\ref{kyr5}) and (\ref{kyr6}), respectively,
\[
a_{ij}^{ +}  = {\frac{{{\sum\limits_{\beta \in J_{r,\,n} {\left\{
{i} \right\}}} {{\left| {\left( {\left( {{\rm {\bf A}}^{ *} {\rm
{\bf A}}} \right)_{\,.\,i} \left({\rm {\bf a}}_{.j}^{ *}\right) }
\right){\kern 1pt} _{\beta }^{\beta} }  \right|}{\kern 1pt}} }
}}{{{\sum\limits_{\beta \in J_{r,\,\,n} } {{\left| {\left( {{\rm
{\bf A}}^{ *} {\rm {\bf A}}} \right){\kern 1pt} {\kern 1pt}
_{\beta} ^{\beta} }  \right|}}} }}}, \quad a_{ij}^{ +}  =
{\frac{{{\sum\limits_{\alpha \in I_{r,m} {\left\{ {j} \right\}}}
{\,{\left| {\left( {\left( {{\rm {\bf A}}{\rm {\bf A}}^{ *} }
\right)_{\,j\,.} \left({\rm {\bf a}}_{i.}^{ *}\right) }
\right)\,_{\alpha} ^{\alpha} } \right|}}} }}{{{\sum\limits_{\alpha
\in I_{r,\,m}} {{\left| {\left( {{\rm {\bf A}}{\rm {\bf A}}^{ *} }
\right){\kern 1pt}  _{\alpha }^{\alpha} }  \right|}}} }}}
\]
\noindent for all $i = 1,\ldots,n $ and $j = 1,\ldots,m $. These
determinantal  representations are original in this case as well.
It is reflected in \cite{ky2}. The analogous result is obtained in
\cite{yu}.
\end{remark}
\begin{remark} The Gram polynomial of a matrix ${\rm {\bf A}}$,
\[
\det({\rm {\bf I}}_{m}+ z{\rm {\bf A}}{\rm {\bf
A}}^{*})=1+a_{1}z+\cdots + a_{m}z^{ m},
\]
is used  for generalizing the Moore-Penrose inverse and the
classical Cramer's rule to the corresponding undetermined and
overdetermined cases over an arbitrary field  in \cite{ge} as
well.
\end{remark}

\section{ Cramer's rule for
a least squares solution of quaternion system linear equations}

\begin{definition}
Suppose
\begin{equation}\label{kyr10}
 {\rm {\bf A}} \cdot {\rm {\bf x}} = {\rm {\bf y}}
 \end{equation}
 is a right
system linear equations over the  quaternion skew field ${\rm
{\mathbb{H}}}$, where ${\rm {\bf A}} \in {\rm
{\mathbb{H}}}^{m\times n}$ is the coefficient matrix, ${\rm {\bf
y}} \in {\rm {\mathbb{H}}}^{m\times 1}$  is a column of constants,
and ${\rm {\bf x}} \in {\rm {\mathbb{H}}}^{n\times 1}$ is a
unknown column.
 The least square solution of (\ref{kyr10}) (with the least norm) is called the vector ${\rm {\bf x}}^{0}$ satisfying
\[
\parallel {{\rm {\bf x}}^{0}}\parallel = {\mathop {\min
}\limits_{{\rm {\bf x}} \in {\rm {\mathbb{H}}}^{n}}} \left\{
\parallel  {{\rm {\bf \tilde {x}}}}\parallel :\,\parallel {{\rm
{\bf A}} \cdot {\rm {\bf \tilde {x}}} - {\rm {\bf y}}}
\parallel  = {\mathop {\min} \limits_{{\rm {\bf x}} \in {\rm
{\mathbb{H}}}^{n}}} \parallel {{\rm {\bf A}} \cdot {\rm {\bf x}} -
{\rm {\bf y}}}\parallel \right\},
\]
where ${\rm {\mathbb{H}}}^{n}$ is  an $n$-dimension right
quaternion vector space.
\end{definition}
We recall that in the right quaternion vector space ${\rm
{\mathbb{H}}}^{n}$ by definition of the inner product of vectors
we put ${\left\langle {{\rm {\bf x}},{\rm {\bf y}}} \right\rangle}
: = {\rm {\bf y}}^{ *} {\rm {\bf x}} = \overline {y_{1}}  \cdot
x_{1}  + \cdots + \overline {y_{n}} \cdot x_{n} $ and ${\left\|
{{\rm {\bf x}}} \right\|}: = \sqrt {{\left\langle {{\rm {\bf
x}},{\rm {\bf x}}} \right\rangle} } $ is the norm
 of a vector ${\rm {\bf x}}\in {\rm {\mathbb{H}}}^{n}$.
By analogy to a complex case (see, e.g. \cite{ga}) we can prove
the following theorem.

\begin{theorem}
The vector ${\rm {\bf x}} = {\rm {\bf A}}^{ +} {\rm {\bf y}}$
   is the least square solution of  (\ref{kyr10}).
\end{theorem}

\begin{definition}
Suppose
\begin{equation}
\label{kyr11} {\rm {\bf x}} \cdot {\rm {\bf A}} = {\rm {\bf y}}
\end{equation} is a left
system linear equations over the  quaternion skew field ${\rm
{\mathbb{H}}}$, where ${\rm {\bf A}} \in {\rm
{\mathbb{H}}}^{m\times n}$ is the coefficient matrix, ${\rm {\bf
y}} \in {\rm {\mathbb{H}}}^{1\times n}$ is a row of constants, and
${\rm {\bf x}} \in {\rm {\mathbb{H}}}^{1\times m}$ is a unknown
row. The least square solution of (\ref{kyr11}) (with the least
norm) is called the vector ${\rm {\bf x}}^{0}$ satisfying
\[{\left\| {{\rm {\bf x}}^{0}}
\right\|} = {\mathop {\min} \limits_{{\rm {\bf \tilde {x}}} \in
{}^{m}{\rm {\mathbb{H}}}}} {\left\{ {{\left\| {{\rm {\bf \tilde
{x}}}} \right\|}:\,\;{\left\| {{\rm {\bf \tilde {x}}} \cdot {\rm
{\bf A}} - {\rm {\bf y}}} \right\|} = {\mathop {\min}
\limits_{{\rm {\bf x}} \in {}^{m}{\rm {\mathbb{H}}}}} {\left\|
{{\rm {\bf x}} \cdot {\rm {\bf A}} - {\rm {\bf y}}} \right\|}}
\right\}},\]
 where ${}^{m}{\rm {\mathbb{H}}}$ is  an $m$-dimension
left quaternion vector space.
\end{definition}

We recall that in the left quaternion vector space ${}^{m}{\rm
{\mathbb{H}}}$ by definition of the inner product of vectors we
can put ${\left\langle {{\rm {\bf x}},{\rm {\bf y}}}
\right\rangle} = {\rm {\bf x}}{\rm {\bf y}}^{ *}  = x_{1} \cdot
\overline {y_{1}} + \cdots + x_{m} \cdot \overline {y_{m}}  $.
Then ${\left\| {{\rm {\bf x}}} \right\|}: = \sqrt {{\left\langle
{{\rm {\bf x}},{\rm {\bf x}}} \right\rangle} }$ is the norm of
${\rm {\bf x}}\in {}^{m}{\rm {\mathbb{H}}}$.
\begin{theorem}\label{kyrc25}
The vector ${\rm {\bf x}} = {\rm {\bf y}} \cdot {\rm {\bf A}}^{ +
}$ is the least square solution of (\ref{kyr11}).
\end{theorem}
\begin{theorem}\label{kyrc18}
\begin{enumerate}
\item[(i)] If $\rank{\rm {\bf A}} = n$, then for the
least square solution  ${\rm {\bf x}}^{0} = (x_{1}^{0} ,\ldots
,x_{n}^{0} )^{T}$ of (\ref{kyr10}), we get for all $j =1,\ldots,n
$
\begin{equation}
\label{kyr12} x_{j}^{0} = {\frac{{{\rm{cdet}} _{j} \left( {{\rm
{\bf A}}^{ *} {\rm {\bf A}}} \right)_{.j} \left( {{\rm {\bf f}}}
\right)}}{{{\rm ddet} {\rm {\bf A}}}}},
\end{equation}
 where ${\rm {\bf f}} = {\rm {\bf A}}^{ *} {\rm {\bf y}}. $
\item[(ii)] If $\rank{\rm {\bf A}} = k \le m < n$, then for all
$j = 1,\ldots,n $ we have
\begin{equation}
\label{kyr13} x_{j}^{0} = {\frac{{{\sum\limits_{\beta \in
J_{r,\,n} {\left\{ {j} \right\}}} {{\rm{cdet}} _{j} \left( {\left(
{{\rm {\bf A}}^{ *} {\rm {\bf A}}} \right)_{\,.\,j} \left( {{\rm
{\bf f}}} \right)} \right){\kern 1pt} {\kern 1pt} _{\beta}
^{\beta} } } }}{{{\sum\limits_{\beta \in J_{r,\,\,n}} {{\left|
{\left( {{\rm {\bf A}}^{ *} {\rm {\bf A}}} \right){\kern 1pt}
{\kern 1pt} _{\beta} ^{\beta} }  \right|}}} }}}.
\end{equation}
\end{enumerate}
\end{theorem}

{\textit{Proof}}. i) If $\rank{\rm {\bf A}} = n$, then ${\rm {\bf
A}}^{ +} $ can be represented by (\ref{kyr8}). Denote ${\rm {\bf
f}}: = {\rm {\bf A}}^{ *} {\rm {\bf y}}$. Representing ${\rm {\bf
A}}^{ +} {\rm {\bf y}}$ by coordinates we obtain (\ref{kyr12}).

\noindent ii) If $\rank{\rm {\bf A}} = k \le m < n$, then by
Theorem  \ref{kyrc16} we represent the matrix ${\rm {\bf A}}^{ +}
$ by (\ref{kyr5}). Representing  ${\rm {\bf A}}^{ +} {\rm {\bf
y}}$ by coordinates we obtain (\ref{kyr13}). $\blacksquare$

\begin{remark}
In a complex case  the following analogs of (\ref{kyr12}) and
(\ref{kyr13}) are obtained respectively in \cite{ky2} for all $j =
1,\ldots,n $,
\[
x_{j}^{0} = {\frac{{\det \left( {{\rm {\bf A}}^{ *} {\rm {\bf A}}}
\right)_{.j} \left( {{\rm {\bf f}}} \right)}}{{\det \left( {{\rm
{\bf A}}^{ *} {\rm {\bf A}}} \right)}}},\; x_{j}^{0} =
{\frac{{{\sum\limits_{\beta \in J_{r,\,n} {\left\{ {j} \right\}}}
{{\left| {\left( {\left( {{\rm {\bf A}}^{ *} {\rm {\bf A}}}
\right)_{.\,j}\, \left( {{\rm {\bf f}}} \right)} \right){\kern
1pt} {\kern 1pt} _{\beta} ^{\beta} }  \right|}}}
}}{{{\sum\limits_{\beta \in J_{r,\,n} } {{\left| {\left( {{\rm
{\bf A}}^{ *} {\rm {\bf A}}} \right) {\kern 1pt} _{\beta} ^{\beta}
}  \right|}}} }}}.\]
\end{remark}
\begin{theorem}\begin{enumerate}
\item[(i)] If $\rank{\rm {\bf A}} = m$, then for  ${\rm {\bf
x}}^{0} = (x_{1}^{0} ,\ldots ,x_{m}^{0} )$ of (\ref{kyr11}) we
obtain for all $i = 1,\ldots,m $

\begin{equation}
\label{kyr14} x_{i}^{0} = {\frac{{{\rm{rdet}} _{i} \left( {{\rm
{\bf A}}{\rm {\bf A}}^{ *} } \right)_{i.} \left( {{\rm {\bf z}}}
\right)}}{{{\rm ddet} {\rm {\bf A}}}}},
\end{equation}
where ${\rm {\bf z}} = {\rm {\bf y}}{\rm {\bf A}}^{ *} $.
\item[(ii)] If
$\rank{\rm {\bf A}} = k \le n < m$, then for all $i = 1,\ldots,m $
we have
\begin{equation}
\label{kyr15} x_{i}^{0} = {\frac{{{\sum\limits_{\alpha \in I_{r,m}
{\left\{ {i} \right\}}} {{\rm{rdet}} _{i} \left( {\left( {{\rm
{\bf A}}{\rm {\bf A}}^{ *} } \right)_{\,i\,.} ({\rm {\bf z}})}
\right)\,_{\alpha} ^{\alpha} } }}}{{{\sum\limits_{\alpha \in
I_{r,\,m}}  {{\left| {\left( {{\rm {\bf A}}{\rm {\bf A}}^{ *} }
\right) {\kern 1pt} _{\alpha} ^{\alpha} } \right|}}} }}}.
\end{equation}
\end{enumerate}
\end{theorem}
 The proof of this theorem is analogous to that
of Theorem \ref{kyrc18}.

\begin{remark}
In a complex case the following analogs of  (\ref{kyr14}) and
(\ref{kyr15}) respectively are obtained  in \cite{ky2}  for all $i
= 1,\ldots,m $,

\[
x_{i}^{0} = {\frac{{\det  \left( {{\rm {\bf A}}{\rm {\bf A}}^{ *}
} \right)\,_{i.} \left( {{\rm {\bf z}}} \right)}}{{\det {\rm {\bf
A}}{\rm {\bf A}}^{ *} }}},\,\,\, x_{i}^{0} =
{\frac{{{\sum\limits_{\alpha \in I_{r,m} {\left\{ {i} \right\}}}
{{\left| {\left( {\left( {{\rm {\bf A}}{\rm {\bf A}}^{ *} }
\right)_{\,i\,.} \left( {{\rm {\bf z}}} \right)}
\right)\,_{\alpha} ^{\alpha} }  \right|}}}
}}{{{\sum\limits_{\alpha \in I_{r,\,m}}  {{\left| {\left( {{\rm
{\bf A}}{\rm {\bf A}}^{ *} } \right){\kern 1pt} {\kern 1pt}
_{\alpha} ^{\alpha} }  \right|}}} }}}.\]
\end{remark}

\section{Example}
Let us consider the left system of linear equations.
\begin{equation}\label{kyr17}
 {\left\{ {{\begin{array}{*{20}c}
 {x_{1} i + 2x_{2} i - x_{3} = i,} \hfill \\
 { - x_{1} k + x_{2} j + x_{3} j = j,} \hfill \\
 {x_{1} j + x_{2} + x_{3} k = k,} \hfill \\
 {x_{1} + x_{2} k + x_{3} i = 1.} \hfill \\
\end{array}} } \right.}
\end{equation}
The coefficient matrix of the system is the matrix ${\rm {\bf A}}
=
\begin{pmatrix}
  i & -k & j & 1 \\
  2i & j & 1 & k \\
  -1 & j & k & i
\end{pmatrix}
$. The row of unknown is ${\rm {\bf x}} = \left(
{{\begin{array}{*{20}c}
 {x_{1}}  \hfill & {x_{2}}  \hfill & {x_{3}}  \hfill \\
\end{array}} } \right)$ and the row of constants is ${\rm {\bf y}} = \left(
{{\begin{array}{*{20}c}
 {i} \hfill & {j} \hfill & {k} \hfill & {1} \hfill \\
\end{array}} } \right)$. Then for (\ref{kyr17}) we have ${\rm {\bf x}} \cdot {\rm
{\bf A}} = {\rm {\bf y}}$. We obtain
\[
{\rm {\bf A}}^{ *}  =\begin{pmatrix}
  -i & -2i & -1 \\
k & -j & -j \\
  -j & 1 & -k \\
  1 & -k & -i
\end{pmatrix}, \]
\[ {\rm {\bf A}}{\rm {\bf A}}^{ *}  =\begin{pmatrix}
4 & 2 - i + j - k & - 4i \\
  2 + i - j + k & 7 & 1 - 2i - j - k \\
  4i & 1 + 2i + j + k & 4
\end{pmatrix}.
\]

Since ${\rm ddet} {\rm {\bf A}} = \det {\rm {\bf A}}{\rm {\bf
A}}^{ *} = {\rm rdet} _{1} {\rm {\bf A}}{\rm {\bf A}}^{ *}  = 0$
and
\[\begin{array}{c}
  \det \left( {{\rm {\bf A}}{\rm {\bf A}}^{ *} } \right)^{33} =
  {\rm
rdet} _{1}\begin{pmatrix}
  4 & 2 - i + j - k \\
  2 + i - j + k & 7
\end{pmatrix} = \\=4 \cdot 7 - \left( {2 - i + j - k} \right) \cdot
\left( {2 + i - j + k} \right) = 21 \ne 0,
\end{array}
\]
 by Lemma
\ref{kyrc6} ${\rm rank}\,{\rm {\bf A}} = 2$. We shall represent
${\rm {\bf A}}^{ +} $ by (\ref{kyr6}).

\[\begin{array}{c}
  {\sum\limits_{\alpha \in I_{2,\,3}}  {{\left| {\left( {{\rm {\bf A}}{\rm
{\bf A}}^{ *} } \right){\kern 1pt} {\kern 1pt} _{\alpha} ^{\alpha}
} \right|}}}  = \det
\begin{pmatrix}
  4 & 2 - i + j - k \\
  2 + i - j + k & 7
\end{pmatrix} + \\
  +\det\begin{pmatrix}
    7 & 1 - 2i - j - k \\
    1 + 2i + j + k & 4 \
  \end{pmatrix}
 + \det\begin{pmatrix}
   4 & - 4i \\
   4i & 4 \
 \end{pmatrix} = 42.
\end{array}
\]

Now we shall calculate $r_{i{\kern 1pt} j} = {\sum\limits_{\alpha
\in I_{2,3} {\left\{ {j} \right\}}} {{\rm rdet} _{j} \left( {({\rm
{\bf A}}{\rm {\bf A}}^{ *} )_{\,j\,.} ({\rm {\bf a}}_{i.}^{ *}  )}
\right)\,_{\alpha} ^{\alpha }} } $ for all $i = \overline {1,4} $
and $j = \overline {1,3} $. To obtain $r_{11} $, we consider the
matrix

\[({\rm {\bf A}}{\rm {\bf A}}^{ * })_{\,1\,.} ({\rm {\bf a}}_{1.}^{
*}  ) =\begin{pmatrix}
  -i & -2i & -1 \\
  2 + i - j + k & 7 & 1 - 2i - j - k \\
  4i & 1 + 2i + j + k & 4
\end{pmatrix}.\]

Then we have
\[\begin{array}{c}
   r_{11} ={\rm rdet} _{1}\begin{pmatrix}
  -i & -2i \\
  2 + i - j + k & 7
\end{pmatrix}+ {\rm rdet}_{1}\begin{pmatrix}
  -i & -1 \\
 4i & 4
\end{pmatrix}=\\
   = - i \cdot 7 - \left( { - 2i} \right) \cdot \left(
{2 + i - j + k} \right) -i\cdot 4-(-1\cdot 4i)= - 2 - 3i - 2j -
2k,
\end{array}
\]
\noindent and so forth. Continuing in the same way, we get
\[
{\rm {\bf A}}^{ +}  = {\frac{{1}}{{42}}}\begin{pmatrix}
 - 2 - 3i - 2j - 2k & 2 - 12i + 2j + 2k &  - 3 + 2i + 2j
- 2k \\
  1 + i + 2j + 6k & - 2 + 2i - 6j - 4k & 1 - i - 6j + 2k \\
  - 2 - i - 6j - k & 6 - 2i + 4j + 2k & - 1 + 2i + j -
6k \\
  6 + i + j + 2k & - 4 + 2i - 2j - 6k & 1 - 6i - 2j + k
\end{pmatrix}.
\]
 We find the least square solution by
means of the matrix method by Theorem \ref{kyrc25}
\[
{\rm {\bf x}}^{0} = {\rm {\bf y}} \cdot {\rm {\bf A}}^{ +}  =
{\frac{{1}}{{42}}}\begin{pmatrix}
  8 + 11i + 3j - 3k,& 12 - 4i - 8j, & 11 - 8i + 3j + 3k
\end{pmatrix}.
\]

Now we shall find the least square solution of (\ref{kyr17}) by
means of Cramer's rule by (\ref{kyr15}). We have $
  {\rm {\bf z}} = {\rm {\bf y}} \cdot {\rm {\bf A}}^{ *}  =\begin{pmatrix}
    2 + 2i, & 3, & 2 - 2i
  \end{pmatrix}.
$ Since
\[({\rm {\bf A}}{\rm {\bf A}}^{ * })_{\,1\,.} ({\rm {\bf z}}) =\begin{pmatrix}
  2 + 2i & 3 & 2 - 2i \\
  2 + i - j + k & 7 & 1 - 2i - j - k \\
  4i & 1 + 2i + j + k & 4
\end{pmatrix},\]
we get
\[
   x_{1}^{0} =\frac{{\rm rdet} _{1}\begin{pmatrix}
    2 + 2i & 3 \\
    2 + i - j + k & 7
  \end{pmatrix} + {\rm rdet} _{1}\begin{pmatrix}
    2 + 2i & 2 - 2i \\
    4i & 4
  \end{pmatrix}}{{\sum\limits_{\alpha \in I_{2,\,3}}  {{\left| {\left( {{\rm {\bf A}}{\rm
{\bf A}}^{ *} } \right){\kern 1pt} {\kern 1pt} _{\alpha} ^{\alpha}
} \right|}}} }
  = \frac{8 + 11i + 3j - 3k}{42}.
 \]

Since $ ({\rm {\bf A}}{\rm {\bf A}}^{ *}) _{\,2\,.} ({\rm {\bf
z}}) =\begin{pmatrix}
4 & 2 - i + j - k & - 4i \\
 2 + 2i & 3 & 2 - 2i \\
  4i & 1 + 2i + j + k & 4
\end{pmatrix}$,
\[ \begin{array}{c}
  x_{2}^{0} = \frac{{\rm rdet} _{2}\begin{pmatrix}
  4 & 2 - i + j - k \\
  2 + 2i & 3
\end{pmatrix} + {\rm rdet} _{1}\begin{pmatrix}
  3 & 2 - 2i \\
  1 + 2i + j + k & 4
\end{pmatrix}}{{\sum\limits_{\alpha \in I_{2,\,3}}  {{\left| {\left( {{\rm {\bf A}}{\rm
{\bf A}}^{ *} } \right){\kern 1pt} {\kern 1pt} _{\alpha} ^{\alpha}
} \right|}}} }
  =\frac{12 - 4i - 8j}{42}.
\end{array}
\]

Since \[ ({\rm {\bf A}}{\rm {\bf A}}^{ *})_{\,3\,.} ({\rm {\bf
z}}) =\begin{pmatrix}
4 & 2 - i + j - k & - 4i \\
  2 + i - j + k & 7 & 1 - 2i - j - k \\
  2 + 2i & 3 & 2 - 2i
\end{pmatrix},\]

 \[\begin{array}{c}
    x_{3}^{0} =\frac{{\rm rdet}_{2}\begin{pmatrix}
   4 & - 4i \\
   2 + 2i & 2 - 2i \
 \end{pmatrix} + {\rm rdet} _{1}\begin{pmatrix}
   7 & 1 - 2i - j - k \\
   3 & 2 - 2i \
 \end{pmatrix}}{{\sum\limits_{\alpha \in I_{2,\,3}}  {{\left| {\left( {{\rm {\bf A}}{\rm
{\bf A}}^{ *} } \right){\kern 1pt} {\kern 1pt} _{\alpha} ^{\alpha}
} \right|}}} }
    =\frac{11 - 8i + 3j + 3k}{42}.\
 \end{array}
\]

\end{document}